\documentclass[preprint,3p,12pt]{elsarticle}
\usepackage{mathrsfs}
\usepackage{amsmath}
\usepackage{stmaryrd}
\usepackage{bbding}
\usepackage{dcolumn}
\usepackage{graphicx}
\usepackage{amsfonts}
\usepackage{amssymb}
\usepackage{psfrag}
\usepackage{wrapfig}
\usepackage{subfigure}
\usepackage{makeidx}
\usepackage{bm}
\usepackage{epsf}
\usepackage{epsfig}
\usepackage{setspace}
\usepackage{graphicx}
\usepackage{epstopdf}
\usepackage{psfrag}
\usepackage{subfigure}
\usepackage{color}

\begin{document}
\title{A High-Order Finite Volume GENO Scheme with Implicit Time Integration for Three-Temperature Radiation Diffusion Equations}

\author[HKUST1]{Fengxiang Zhao}
\ead{fzhaoac@connect.ust.hk}

\author[HKUST1]{Yaqing Yang}
\ead{yangyq@ust.hk}

\author[IAPCM]{Yibing Chen}
\ead{chen\_yibing@iapcm.ac.cn}

\author[HKUST1,HKUST2]{Kun Xu\corref{cor}}
\ead{makxu@ust.hk}

\address[HKUST1]{Department of Mathematics, Hong Kong University of Science and Technology, Clear Water Bay, Kowloon, Hong Kong, China}
\address[IAPCM]{Institute of Applied Physics and Computational Mathematics and National Key laboratory of Computational physics, and Center for Applied Physics and Technology, Peking University, Beijing, China}
\address[HKUST2]{Shenzhen Research Institute, Hong Kong University of Science and Technology, Shenzhen, China}
\cortext[cor]{Corresponding author}

\begin{abstract}
This study presents a high-order finite volume scheme capable of large time-step integration for three-temperature radiation diffusion (3TRD) equations, where conservation is naturally achieved through energy update. To handle local large gradients and discontinuities in temperature, a central generalized ENO (GENO) reconstruction is developed for diffusion systems, which achieves essentially non-oscillatory reconstruction for discontinuous solutions.
Compared to conventional nonlinear reconstruction methods, its most distinctive feature is the central-type symmetric sub-stencils, which ensure consistency between the numerics and the isotropic nature of thermal diffusion.
Additionally, the central GENO method provides smooth states of temperature and temperature gradient at interfaces, facilitating the evaluation of numerical fluxes. Furthermore, interface flux evaluation for cases with discontinuous physical property parameters is modeled.
To address the extremely small time-step issue caused by stiff diffusion and source terms, a dual-time-stepping method based on implicit time discretization is developed for the first time in 3TRD systems, with the advantage of decoupling temporal discretization from complex nonlinear spatial discretization.
A series of numerical examples validates the high accuracy, physical property preservation, strong robustness, and large time-step integration capability of the present high-order central GENO scheme.
\end{abstract}


\maketitle

\section{Introduction}

The three-temperature radiative diffusion (3TRD) equations \cite{zeyao2004parallel,jiang2007some} constitute a fundamental model in inertial confinement fusion (ICF) and high-energy-density physics for describing the radiative energy transport in multi-material systems and the energy exchange among electrons, ions, and photons.
The 3TRD equations are characterized by severe nonlinearity and multi-scale stiffness. The nonlinearity stems from the pronounced temperature dependence of the diffusion and energy-exchange coefficients, most notably the quartic or higher-order power-law scaling of the radiation diffusion coefficient. Meanwhile, the multi-scale stiffness arises from the vast disparity in characteristic time scales among various physical processes. For instance, the speed of radiative diffusion can outpace that of material thermal conduction by several orders of magnitude.
Moreover, realistic scenarios typically involve multi-material systems in which material parameters (such as thermal conductivity and specific heat capacity) are discontinuous at material interfaces, causing temperature profiles to exhibit extremely steep gradients across these interfaces. In addition, numerical studies must ensure fundamental physical constraints, including bound preservation and positivity of temperatures, along with energy conservation.
These inherent properties present substantial challenges for numerical simulations with respect to achieving high accuracy, computational efficiency, and algorithmic robustness.

Finite volume methods with fully implicit time integration constitute the dominant strategy for solving 3TRD equations \cite{yu2019finite,peng2020positivity}. Existing numerical schemes are restricted to spatial discretizations of at most second-order accuracy. While high-order methods offer substantial advantages in computational efficiency and solution accuracy \cite{wang2007high}, their extension to 3TRD systems remains largely unexplored.
Consequently, the development of robust and efficient high-order schemes capable of simultaneously preserving physical constraints and handling strong nonlinearities represents a significant open challenge. This work addresses this gap by proposing a high-order implicit scheme specifically tailored for the 3TRD system.

Nonlinear spatial reconstruction is fundamental to constructing high-order schemes, critically influencing accuracy and stability for problems involving multi-scale spatial distributions and extreme conditions such as severely steep temperature gradients and near-zero temperatures. The objective of nonlinear reconstruction is to adaptively transition from high-order reconstruction in smooth regions to robust lower-order reconstruction near large gradients and discontinuities. WENO reconstruction, a widely adopted nonlinear method in compressible flow simulations, achieves this through nonlinear combinations of lower-order candidate polynomials, adaptively recovering high-order linear reconstruction in smooth regions while reducing to upwind-biased lower-order reconstruction that avoids spurious oscillations near discontinuities \cite{liu,jiang}. Optimization efforts for WENO methods have primarily focused on designing nonlinear weights and constructing optimal candidate polynomials \cite{WENO-Z,zhu2016new}, as these components critically determine performance of WENO schemes.
GENO reconstruction method was originally developed for high-order schemes in compressible flow simulations \cite{zhao2025generalized}. Unlike WENO methods, GENO employs a path function to directly connect high-order linear reconstruction with robust lower-order reconstruction (e.g., second-order ENO or TVD), maintaining linear reconstruction in high-wavenumber regions while adaptively transitioning to lower-order reconstruction at discontinuities. This method exhibits low sensitivity to lower-order candidate polynomials, making it particularly well-suited for spatially multi-scale problem and efficient nonlinear reconstruction.

This study presents a novel central GENO reconstruction for 3TRD equations. The central GENO is distinguished by its construction at cell interfaces of both a high-order linear polynomial and a second-order polynomial, based on the central-type large stencil and the sub-stencil, respectively, thereby yielding single high-order state of temperature and its derivative at interfaces. The central-type reconstruction ensures compatibility of the spatial discretization with isotropic diffusion physics, while the smooth interface reconstruction facilitates evaluation and modeling of numerical fluxes. The GENO path function provides a methodological framework for achieving high-order nonlinear reconstruction using only one high-order linear reconstruction and one second-order reconstruction, thereby circumventing the cumbersome construction of multiple sub-stencil polynomials.
Theoretically, the central GENO reconstruction possesses essentially non-oscillatory properties and, in numerical experiments, demonstrates bound-preserving characteristics for discontinuous problems. Moreover, to facilitate simple and flexible treatment of 3D problems, this study proposes a ``1D $+$ 2D'' combined reconstruction approach. The interface average is first obtained via 1D reconstruction along the interface normal direction. For linear diffusion problems, high-order spatial discretization can be achieved without requiring reconstruction in other directions. For nonlinear problems, 2D reconstruction is subsequently performed, using the interface average from the 1D reconstruction, to determine values at Gaussian quadrature points on the interface.

To overcome the problem of extremely small time steps caused by strong stiffness, implicit discretization schemes must be employed for the time integration of 3TRD, such as the first-order backward Euler scheme or the second-order Crank-Nicolson scheme. Implicit discretization yields a system of nonlinear equations, which typically requires further linearization treatment. The solution of this system can be achieved using methods such as Picard iteration \cite{yu2019finite}, Krylov subspace iteration methods \cite{saad1986gmres,brown1990_GMRES}, incomplete factorization methods, and algebraic multigrid (AMG) methods. To achieve rapid convergence and efficient computation, these methods are often used in combination, for example, Krylov subspace iteration methods preconditioned by either ILU methods or AMG methods.

This study employs the dual time-stepping method to transform and solve the implicit discretization equations. The key technique is to introduce a pseudo-time derivative term, thereby converting the implicit discretization problem in physical time into a pseudo-temporal evolution problem toward steady-state convergence for time-independent variables, where the corresponding converged solution represents the unsteady solution at that physical time step \cite{jameson1991time,pulliam1993time,tan2017time}. The dual time-stepping method has been widely applied to solve unsteady flow problems. High-order backward difference schemes can be employed in physical time to ensure temporal accuracy, while mature implicit steady-state solution techniques are utilized in the pseudo-time direction to iteratively solve the discretization equations until the pseudo-time residuals converge sufficiently. The flux evaluation in this study is based on high-order nonlinear reconstruction method. The complexity of high-order reconstruction makes it impossible to explicitly formulate the implicit discretization equations. Therefore, within the conventional implicit time discretization framework for 3TRD, it is difficult to implement high-order spatial discretization. The application of the dual time-stepping method to the 3TRD equations to develop a new implicit scheme represents one of the core contributions of this study. The benefit brought by the dual time-stepping method is the decoupling of implicit time discretization and nonlinear spatial discretization.

This paper is organized as follows.
Section 2 presents the finite volume method for 3TRD equations.
The high-order GENO reconstruction will be introduced in Section 3.
Section 4 presents the modeling and calculation of fluxes as well as the discretization of source terms.
Section 5 presents the implicit acceleration method.
Section 6 provides validation test cases, and Section 7 concludes the paper.

\section{3TRD equations and finite volume method}

The 3TRD equations \cite{zeyao2004parallel,yu2019finite} are given as
\begin{equation}\label{3T-pde}
\frac{\partial \textbf{W}}{\partial t} -\nabla\cdot \textbf{F} = \textbf{S},
\end{equation}
where $\mathbf{W}$ denotes the energy vector of electrons, ions, and photons, and $\mathbf{F}$ is the corresponding flux.
The source term $\mathbf{S}$ accounts for the energy transfer among electrons, ions, and photons.
The specific forms of $\mathbf{W}$, $\mathbf{F}$, and $\mathbf{S}$ are
\begin{equation*}
{\textbf{W}} =
\left(
\begin{array}{c}
c_eT_e\\
c_iT_i\\
c_rT_r\\
\end{array}
\right), \\
{\textbf{F}} =
\left(
\begin{array}{c}
k_e\nabla T_e\\
k_i\nabla T_i\\
k_r\nabla T_r\\
\end{array}
\right),\\
\end{equation*}
and
\begin{equation*}
{\textbf{S}} =
\left(
\begin{array}{c}
\omega_i(T_i-T_e)+\omega_r(T_r-T_e)\\
\omega_i(T_e-T_i)\\
\omega_r(T_e-T_r)\\
\end{array}
\right).
\end{equation*}
where $k_{\alpha}$ and $c_{\alpha}$ $(\alpha\in \{e,i,r\})$ represent the diffusion coefficient and the volumetric heat capacity of species $\alpha$, respectively.
For electrons and ions, the volumetric heat capacity is defined by \(c_{\alpha}=\rho\,c_{v\alpha}\) (\(\alpha\in\{e,i\}\)), where \(c_{v\alpha}\) is the specific heat capacity. For photons, there is \(c_r = 4c_{vr}T_r^{3}\). The parameter \(\omega_{\alpha}\) denotes the inter-species heat-exchange coefficient.

Integrating Eq. (\ref{3T-pde}) over a control volume to construct the finite volume scheme yields
\begin{equation}\label{3T-fvm}
\frac{\text{d}\textbf{W}_{j}}{\text{d}t}=\frac{1}{\big|\Omega_j\big|}\int_{\partial \Omega_j} \textbf{F}\cdot \textbf{n} \mathrm{d} \Gamma +\int_{\Omega_j} \mathbf{S} \mathrm{d} \Omega,
\end{equation}
where $\textbf{W}_{j}$ are the cell-averaged energy variables, $\textbf{F}$ are the numerical fluxes at cell interfaces, $\big|\Omega_j\big|$ is the volume of $\Omega_j$, and $\textbf{n}$ is the unit outer normal vector to the interface $\partial \Omega_j$.
The cell-averaged variables $\textbf{W}_{j}$ are defined as
\begin{align*}
\textbf{W}_{j}= \frac{1}{\big| \Omega_j \big|} \int_{\Omega_j} \textbf{W}(\mathbf{x}) \mathrm{d}\Omega.
\end{align*}
The integral on the cell interfaces on the right-hand side of Eq. (\ref{3T-fvm}), which is a line integral in two dimensions and a surface integral in three dimensions, is discretized using Gaussian quadrature as
\begin{align*}
\int_{\partial \Omega_j} \textbf{F}\cdot \textbf{n} \mathrm{d} \Gamma=\sum_{l=1}^{l_0}\big( \big|\Gamma_{l} \big| \sum _{k=1}^{k_0} w_k \textbf{F}(\mathbf{x}_k)\cdot \textbf{n}_l \big).
\end{align*}
Here, \(l_0\) denotes the number of interfaces of cell \(\Omega_j\), and \(|\Gamma_l|\) is the length (in 2D) or
area (in 3D) of the \(l\)-th interface. Moreover, \(k_0\) and \(w_k\) are the number of Gaussian quadrature
points and the quadrature weights, respectively.

For the 3TRD problem, the development of high-order finite volume schemes requires the reconstruction of temperature states at interfaces to compute numerical fluxes, the cell integration of source terms, and efficient time integration. The high-order spatial reconstruction and efficient time integration constitute the primary focus of this work and will be elaborated upon in the following sections.

\section{High-order spatial reconstruction}

High-order spatial reconstruction for temperature is presented in this section.
To achieve non-oscillatory, robust computations near temperature discontinuities and steep gradients, nonlinear reconstruction methods are necessary. This study introduces central GENO reconstruction, applying GENO to diffusion problems for the first time.
For 3D structured grids, a novel ``1D + 2D'' stage-by-stage reconstruction strategy is proposed, featuring high-order accuracy and algorithmic simplicity, making it well-suited for high-order finite volume schemes.
For nonlinear system, the second-stage 2D reconstruction is necessary to achieve high-order convergence.

\subsection{1D Central GENO reconstruction}

\begin{figure}[!htb]
\centering
\hspace*{-5em}
\includegraphics[width=0.70\textwidth]{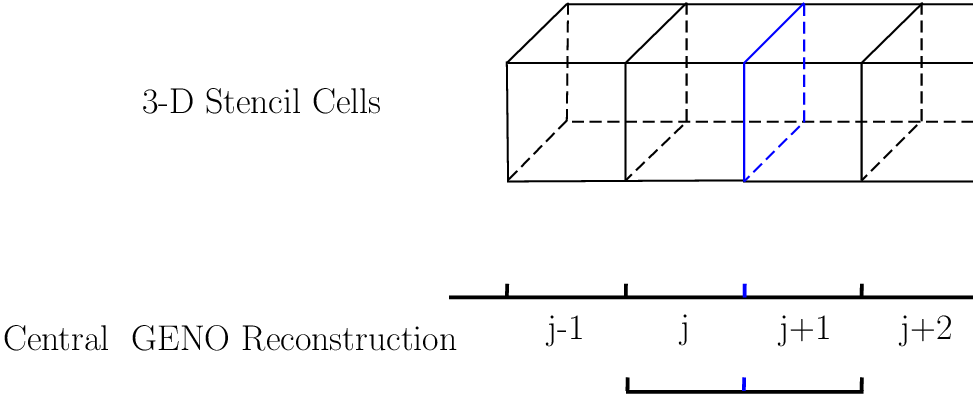}
\caption{\label{1-stencil-1d} Schematic of high-order reconstruction on 3D structured grids. The stencil cells used in the 1D central GENO reconstruction for the ``1D + 2D'' reconstruction strategy. The face with blue edges denotes the target interface to be reconstructed.}
\end{figure}

The 3D high-order reconstruction is performed in two stages. In the first stage, a 1D reconstruction is carried out at cell interfaces to obtain the temperature state on the interface. Figure \ref{1-stencil-1d} illustrates the stencil employed for the 1D reconstruction, including the sub-stencil and the large stencil utilized in the novel central GENO method.
The linear reconstruction that underpins the nonlinear method is presented first.
Based on the symmetric large stencil, a linear fourth-order polynomial is constructed, with the reconstructed value and its derivative at the interface given by
\begin{equation}\label{GENO-linear-p3}
\begin{split}
&p^3(x_{j+1/2})=(-Q_{j-1}+7Q_{j}+7Q_{j+1}-Q_{j+2})/12,\\
&p^3_x(x_{j+1/2})=(Q_{j-1}-15Q_{j}+15Q_{j+1}-Q_{j+2})/12h,
\end{split}
\end{equation}
where $Q$ denotes any scalar variable to be reconstructed and $h$ represents the mesh spacing.

In the GENO-based finite volume scheme for the 3TRD equations, only a single two-cell central sub-stencil is employed, as illustrated in Figure \ref{1-stencil-1d}. Notably, this central stencil is absent from conventional ENO or WENO formulations. Its adoption in this work is motivated by its ability to consistently yield physically admissible reconstructions while maintaining natural compatibility with isotropic thermal diffusion.
The resulting second-order interface reconstruction based on this central two-cell sub-stencil is given by
\begin{equation}\label{GENO-linear-p1}
\begin{split}
&p^1(x_{j+1/2})=(Q_{j}+Q_{j+1})/2,\\
&p^1_x(x_{j+1/2})=(Q_{j+1}-Q_{j})/h.
\end{split}
\end{equation}
If a discontinuity (or a steep gradient) appears at any cell interface within the associated large stencil, namely at \(x_{j-1/2}\), \(x_{j+1/2}\), or \(x_{j+3/2}\), the GENO method adaptively reverts to the second-order reconstruction Eq. \eqref{GENO-linear-p1} at \(x_{j+1/2}\).

The GENO method is a nonlinear high-order reconstruction that connects a high-order linear reconstruction with a reliable or physically consistent lower-order reconstruction (e.g., ENO and TVD). This method enables adaptive high-order accuracy for smooth solutions while recovering the lower-order reconstruction when necessary to maintain robustness.
For a comprehensive analysis of the GENO method, readers are referred to \cite{zhao2025generalized}. Based on the GENO formulation, a novel central GENO method is developed for the present diffusion problem by combining the fourth-order polynomial $p^3$ with the second-order polynomial $p^1$, expressed as
\begin{equation}\label{GENO-nonlinear}
\begin{split}
R(\mathbf{x})=&\chi p^3(\mathbf{x}) +(1-\chi) p^1(\mathbf{x}),\\
\vspace{15pt}
\chi=&\mathrm{\mathbf{Tanh}}(C \alpha)/ \mathrm{\mathbf{Tanh}}(C),
\end{split}
\end{equation}
where $\chi$ is the path function and the parameter $C$ is fixed at $20$.
The ultimate smoothness indicator $\alpha$, which quantifies the smoothness of the linear high-order reconstruction, is given as
\begin{equation}\label{GENO-IS-linear}
\begin{split}
\alpha&=\frac{2\alpha^H}{\alpha^H+\alpha^L},\\
\vspace{10pt}
\alpha^H&=1+\big(\frac{IS^{\tau}}{IS^H+\epsilon}\big)^r, \alpha^L=1+\big(\frac{IS^{\tau}}{IS^L+\epsilon}\big)^r,
\end{split}
\end{equation}
where the power $r$ is $2$, and the small parameter $\epsilon$ is taken as $10^{-15}$.
Building upon the non-compact GENO reconstruction in \cite{zhao2025generalized} and taking into account the symmetric nature of the proposed central GENO method, $IS^H$ and $IS^L$ are defined by the smoothness indicators associated with four auxiliary sub-stencils. $IS^{\tau}$ is a smoothness metric associated with higher-order derivative terms of the reconstruction over the large stencil.
The four sub-stencils used to determine $IS^L$ and $IS^H$ are
\begin{equation*}
\begin{split}
&s^*_1=\{\Omega_{j-1},\Omega_{j}\},~~~~~~~~~~s^*_2=\{\Omega_{j+1},\Omega_{j+2}\},\\
&s^*_3=\{\Omega_{j-1},\Omega_{j},\Omega_{j+1}\},~~~s^*_4=\{\Omega_{j},\Omega_{j+1},\Omega_{j+2}\}.
\end{split}
\end{equation*}
The smoothness indicators $IS^*_k$ on $s^*_k$ are calculated following the conventional definition \cite{jiang} as
\begin{equation*}
\begin{split}
&IS^*_1=(Q_{j}-Q_{j-1})^2, \\
&IS^*_2=(Q_{j+2}-Q_{j+1})^2,\\
&IS^*_3=\frac{13}{12}(Q_{j-1}-2Q_{j}+Q_{j+1})^2+\frac{1}{4}(Q_{j-1}-Q_{j+1})^2,\\
&IS^*_4=\frac{13}{12}(Q_{j+2}-2Q_{j+1}+Q_{j})^2+\frac{1}{4}(Q_{j+2}-Q_{j})^2.
\end{split}
\end{equation*}
Finally, the resulting $IS^L$, $IS^H$ and $IS^{\tau}$ for the GENO method are given by
\begin{equation}\label{GENO-1d-IS-HL}
\begin{split}
&IS^L=\mathrm{\mathbf{Min}}\{IS^*_k \mid k=1,2,3,4\},~~IS^H=\mathrm{\mathbf{Max}}\{IS^*_k \mid k=1,2,3,4\},\\
&IS^{\tau}=|IS^*_3-IS^*_4|.
\end{split}
\end{equation}

\subsection{2D GENO reconstruction at quadrature points}

\begin{figure}[!htb]
\centering
\includegraphics[width=0.35\textwidth]{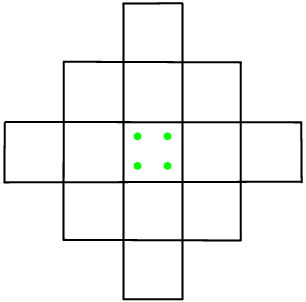}
\caption{\label{1-stencil-2d} Schematic of high-order reconstruction on 3D structured grids. 2D reconstruction stencil for the ``1D + 2D'' stage-by-stage strategy. Green dots denote Gaussian quadrature points on interfaces for flux evaluation.}
\end{figure}

The second stage of the novel ``1D + 2D'' dimensionally split strategy employs a 2D reconstruction to provide the temperature states at quadrature points on the cell interface for fluxes evaluation, using the stencil shown in Figure \ref{1-stencil-2d}. A cubic polynomial $q^3(y,z)$ is determined via least-squares method subject to the following reconstruction constraints
\begin{align}\label{Recons-2-stage}
\big(\frac{1}{\big|\Omega_m \big|} \int_{ \Omega_m} \varphi_k(\bm{x}) \mathrm{d}x \mathrm{d}y \big) a_k=Q_m, ~m=0,1,\cdots,12,
\end{align}
where summation over the repeated index $k$ is implied (Einstein summation convention), $Q$ is any reconstructed scalar component, and $\big|\Omega_m \big|$ is the cell area. For grids with unequal spacing in the two directions, reconstruction is performed on a transformed uniform computational grid.

In multidimensional reconstruction, exact equality between the number of coefficients $a_k$ and reconstruction constraints is generally unattainable. A constrained least-squares approach is employed where the cell average is strictly preserved by $q^3$, while remaining conditions are satisfied in the least-squares sense, yielding
\begin{equation} \label{CLS-system}
\left(
\begin{array}{cc}
\mathbf{A}_{0,k} & 0 \\
2\mathbf{A}_{m,k} \mathbf{A}_{m,n} & \mathbf{A}_{0,n}^\mathrm{T} \\
\end{array}
\right)
\left(
\begin{array}{c}
\mathbf{a}_k \\
c \\
\end{array}
\right)
=
\left(
\begin{array}{cc}
1 & \mathbf{0} \\
\mathbf{0} & 2\mathbf{A}_{m,n}^\mathrm{T} \\
\end{array}
\right)
\mathbf{b},
\end{equation}
where $m=1,2,\ldots,12$ and $k,n=1,2,\ldots,10$. The matrix $\mathbf{A}$ is formed from the cell averages of the basis functions given in Eq. (\ref{Recons-2-stage}), resulting in a $13\times10$ matrix. The vector $\mathbf{b}$ is a $13\times1$ column vector whose components are the cell-averaged values $Q_m$. The parameter $c$ represents an auxiliary Lagrange multiplier introduced to enforce the cell-average constraint.

To handle discontinuities, the GENO method is similarly employed. The lower-order component uses a simple, computationally efficient second-order ENO reconstruction. The four three-cell sub-stencils are:
\begin{equation*}
\begin{split}
s_1=\{\Omega_0,\Omega_1,\Omega_2\},~~~s_2=\{\Omega_0,\Omega_2,\Omega_3\},\\
s_3=\{\Omega_0,\Omega_3,\Omega_4\},~~~s_4=\{\Omega_0,\Omega_4,\Omega_1\},
\end{split}
\end{equation*}
where $\Omega_k$ ($k=1,2,3,4$) denotes the face-neighboring cells of $\Omega_0$ arranged in clockwise or counterclockwise order.
In the present tangential 2D reconstruction, biased sub-stencils are employed. Numerical experiments indicate that maintaining isotropy in the tangential nonlinear reconstruction is not essential.
The 2D GENO reconstruction follows the multi-dimensional formulation described in \cite{zhao2025generalized}.
For the sake of completeness, the formulas used for computation are briefly summarized below.
\begin{equation}\label{GENO-2d-IS}
\begin{split}
&IS^{L}=\big(\sum_{k=1}^{4}IS_k -\mathrm{\mathbf{Max}}\{IS_1,IS_3\} -\mathrm{\mathbf{Max}}\{IS_2,IS_4\} \big)/2,\\
&IS^{H}=IS^{q^3}= \sum_{k=2}^{10}a_k^2, ~~\widetilde{IS}^{H}=IS^{q^3,*} = \sum_{k=2}^{6}a_k^2,\\
&IS^{\tau}=|IS^{H}-\widetilde{IS}^{H}|.
\end{split}
\end{equation}
Here, $\widetilde{IS}^{H}$ differs from $IS^{H}$ by the omission of terms involving the third-order derivatives of the polynomial.
By employing the Taylor basis $\varphi_k$ in Eq. (\ref{Recons-2-stage}), both $IS^H$ and $\widetilde{IS}^{H}$ reduce to concise forms determined explicitly by the polynomial coefficients. Furthermore, in the context of multidimensional reconstruction, the parameter $r=3$ is typically adopted to evaluate the GENO path function via Eq. (\ref{GENO-IS-linear}).

\subsection{High-accuracy boundary treatment}

In 3TRD problems, steep temperature gradients are prevalent at the boundaries of the computational domain, necessitating accurate boundary reconstruction. To address this, we implement high-order reconstruction methods that strictly preserve physical constraints. Since the diffusion flux is aligned with the boundary-normal direction, and to maintain algorithmic simplicity, our high-order treatment modifies the temperature profile solely along the normal direction. Consequently, boundary values are evaluated using a purely 1D high-order reconstruction along the normal direction.
Specifically, we detail the reconstruction procedures for two distinct types of boundary conditions below.

Under Dirichlet boundary conditions, the boundary temperature is prescribed as \(Q_b\). Let \(I_{j+1}\) and \(I_{j+2}\) in Figure \ref{1-stencil-1d} denote the computational cells immediately adjacent to the left boundary. Consequently, the stencil adopted for the 1D boundary reconstruction is \(\{Q_b, Q_1, Q_2\}\). To strictly enforce this boundary condition, we propose the following reconstruction conditions:
\begin{align*}
\begin{split}
&\big(\frac{1}{\big|I_{m}\big|} \int_{I_{m}} \varphi_k(\bm{x}) \mathrm{d}x \big) a_k=Q_m, ~m=1,2, \\
&\varphi_k(x_{1/2}) a_k=Q_b,
\end{split}
\end{align*}
where the repeated index \(k\) in \(\varphi_k\) and \(a_k\) denotes the Einstein summation convention.
The normal temperature gradient is readily obtained as
\begin{align*}
\partial Q_{1/2}/\partial x=-\big(6Q_b-7Q_{1}+Q_{2} \big)/2\Delta x.
\end{align*}
For the right boundary case, by symmetry, it is readily obtained as $\partial Q_{j+1/2}/\partial x= -\partial Q_{1/2}/\partial x$. Numerical tests demonstrate that this linear high-order treatment can simultaneously achieve bound-preserving properties.

Under Neumann boundary conditions, the temperature gradient at the boundary is prescribed as \(Q_{b,x}\). Following a similar approach for the left boundary, the 1D reconstruction stencil is defined as \(\{Q_{b,x}, Q_1, Q_2\}\). By constraining the reconstruction polynomial to strictly satisfy these three values, the boundary temperature is readily obtained as:
\begin{align*}
Q_{1/2}=\big(7Q_{1}-Q_{2}-2\Delta x Q_{b,x} \big)/6.
\end{align*}
By symmetry, the temperature at the right boundary is given by:
\begin{align*}
Q_{j+1/2}=\big(7Q_{j}-Q_{j-1}+2\Delta x Q_{b,x} \big)/6.
\end{align*}
These interface temperatures are subsequently utilized to evaluate the temperature-dependent heat conduction and energy-exchange coefficients in practical applications.
In addition, the temperature gradient for diffusion fluxes takes the value prescribed by Neumann boundary conditions.

\section{Flux modeling and source terms evaluation}

This section presents the flux calculation and source term discretization for the 3TRD equations. The flux calculation is based on the GENO reconstruction described in the previous section. We first analyze the treatment of discontinuities using the novel central GENO reconstruction. The GENO reconstruction provides smooth interface temperature states, which facilitates flux calculation; however, the discontinuous material properties at the interfaces require special treatment.

\subsection{ENO property of central GENO and flux modeling}

A novel central GENO method is proposed for handling discontinuities and steep gradients, which is characterized by the dominant central second-order reconstruction from the interface-symmetric sub-stencil for discontinuities. The properties of the central GENO reconstruction are first presented here.
Additionally, in realistic problem simulations, material properties in flux evaluations are temperature-dependent. Consequently, accurate reconstruction of both temperature variables and temperature gradients is essential.

In smooth regions, GENO achieves high-order accurate approximation of the exact temperature field, yielding
\begin{equation*}
R(x)=Q(x)+O(h^4).
\end{equation*}
At temperature discontinuities located at any interface within the large stencil of Figure \ref{1-stencil-1d}, GENO adaptively produces temperature and its gradient via the central 2nd-order reconstruction on the symmetric sub-stencil, giving
\begin{equation*}
\begin{split}
&R(x_{j+1/2})=(Q_{j}+Q_{j+1})/2, \\
&R_x(x_{j+1/2})=(Q_{j+1}-Q_{j})/h.
\end{split}
\end{equation*}
For smooth solutions, the total variation of the high-order linear reconstruction satisfies
\begin{equation}\label{Analysys-GENO-TVD-1}
\begin{split}
TVD\big(\sum_jR_j(x)\big)&=TVD\big(\sum_jQ_j(x)+O(h^4)\big) \\
&=TVD\big(Q(x)\big)+O(h^4).
\end{split}
\end{equation}
While for a fully second-order reconstruction, the total variation satisfies
\begin{equation}\label{Analysys-GENO-TVD-2}
\begin{split}
TVD\big(\sum_j R_{j+1/2}(x)\big)&=\sum_j |Q_{j+1}-Q_{j}| \\
                                         &=\sum_j |Q(x_{j+1})-Q(x_{j})|+O(h^2) \\
                                         &=TVD\big(Q(x)\big)+O(h^2),
\end{split}
\end{equation}
provided that boundary conditions of at least second-order accuracy are imposed.
Thus, the central GENO method achieves ENO properties when handling discontinuities. For general cases involving both discontinuities and smooth regions, the ENO property remains valid because the GENO reconstruction is obtained through a convex combination of the high-order linear reconstruction and the second-order reconstruction via the path function.
The second stage of the proposed ``1D+2D'' splitting reconstruction employs cell-centered GENO reconstruction. Its ENO property has been established in previous studies \cite{zhao2025generalized} and is not detailed here.

Furthermore, in the vicinity of temperature discontinuities, the GENO method provides a locally bounded temperature state at the cell interfaces, along with a physically consistent temperature gradient,
\begin{equation}\label{Analysys-GENO-preserving}
\begin{split}
R(x_{j+1/2})&\in [Q_{j},Q_{j+1}],\\
\mathrm{\mathbf{Sign}}\big(R_x(x_{j+1/2})\big)&=\mathrm{\mathbf{Sign}}\big(Q_{j+1}-Q_{j}\big),
\end{split}
\end{equation}
where the second equation ensures the correct heat flux direction.

The central GENO reconstruction produces temperature states that are smooth at interfaces, enabling direct flux evaluation from its definition as
\begin{equation*}
\widehat{\mathbf{F}}\equiv \mathbf{F}\cdot \mathbf{n}=\hat{k}\frac{\partial \mathbf{T}}{\partial \mathbf{n}}.
\end{equation*}
This indicates an advantage of the central GENO method in that it avoids the need for special treatment of temperature discontinuities at interfaces during flux evaluation.

However, for the extreme case involving material property discontinuities at material interfaces, a specific modeling approach is required.
At a general material interface, $k_{l,r}$ and $\mathbf{T}{l,r}$ denote the thermal conductivities and temperatures in the cells on the left and right sides of the interface, respectively. Based on the heat flux continuity condition, we have
\begin{equation*}
\widehat{\mathbf{F}}=\hat{k}\frac{\partial \mathbf{T}}{\partial \mathbf{n}}=k_l\frac{\mathbf{T}-\mathbf{T}_l}{h/2}=k_r\frac{\mathbf{T}_r-\mathbf{T}}{h/2},
\end{equation*}
where $\mathbf{T}$ denotes the exact interface temperature. At material interfaces, the central GENO reconstruction effectively models continuous temperature and temperature gradient states despite the underlying gradient discontinuities. The equivalent thermal conductivity is obtained as
\begin{equation}
\hat{k}=\frac{2k_l k_r}{k_l+ k_r}.
\end{equation}

\subsection{Source terms evaluation}

The source term in the radiative diffusion equations can be written in a general form as
\begin{equation*}
S=\omega(Q)Q,
\end{equation*}
which is discretized via cell integration as
\begin{equation}\label{S-space-method}
\begin{split}
S_j&=\frac{1}{|\Omega_j|}\int_{\Omega_j}\omega(Q)Q \mathrm{d}\Omega \\
&=\frac{1}{|\Omega_j|} \int_{\Omega_j}\big( \omega(Q_0)+\delta \mathbf{x}\cdot \nabla \omega(Q_0) +O(h^2) \big)\big( Q_0 +\delta \mathbf{x}\cdot \nabla Q_0 +O(h^2)  \big) \mathrm{d}\Omega\\
&=\omega(Q_0) Q_0 +\sum_m \frac{1}{12}\delta x_m^2 \partial_m \omega(Q_0) \partial_m Q_0 +O(h^2).
\end{split}
\end{equation}
The above second-order discretization is exact for constant $\omega$ and for cases where both $\omega$ and $Q$ are linear on each cell.
While fourth-order accuracy can be achieved via Gaussian quadrature, this is not employed here due to considerations of computational efficiency and implementation simplicity.
The temperature gradient on a cell is determined from the gradients at both interfaces through a nonlinear limiter, given by
\begin{equation*}
\nabla Q_0=\mathrm{S_L}\big(\nabla Q_{-1/2},\nabla Q_{1/2}\big),
\end{equation*}
where $\mathrm{S_L}$ denotes the limiter function, and the minmod limiter is employed in this study.

\section{Implicit method for temporal integration}

In this section, we propose a new implicit acceleration method to alleviate the strict time-step constraints caused by stiff source and heat convection terms.
Unlike conventional Picard iteration methods \cite{yu2019finite} for solving implicit discrete schemes, we employ the dual time-stepping method to handle implicit schemes in this study. This approach facilitates obtaining implicit solutions when nonlinear high-order reconstruction is employed in space.
For Eq. (\ref{3T-fvm}), the fully discrete scheme using one-step backward Euler temporal discretization is written as
\begin{equation}\label{implicit-1}
\frac{\mathbf{W}^{n+1}_j-\mathbf{W}^{n}_j}{\Delta t}-\mathcal{L}_j^{n+1}=0,
\end{equation}
where $\mathcal{L}_j$ represents the flux and source terms, given by
\begin{equation}\label{3TRD-L}
\mathcal{L}_j=\frac{1}{\big|\Omega_{j} \big|} \sum_{l=1}^{l_0}\big( \big|\Gamma_{l} \big| \sum _{k=1}^{k_0} w_k \textbf{F}(\mathbf{x}_k)\cdot \textbf{n}_l \big)-\mathbf{S}_j.
\end{equation}

In the large time-step advancement from $t^n$ to $t^{n+1}$, to obtain the high-accuracy unsteady solution of Eq. (\ref{implicit-1}), the dual time-stepping method transforms the problem of solving nonlinear equations into a convergence problem for the time-independent variable $\Delta \bm{W}^m_j = \bm{W}^{m+1}_j - \bm{W}^m_j$ by introducing a pseudo-time derivative term, as follows:
\begin{equation}\label{implicit-2}
\frac{\mathbf{W}^{m+1}_j-\mathbf{W}^{n}_j}{\Delta t}-\mathcal{L}^{m+1}_j=\frac{\mathbf{W}^{m}_j-\mathbf{W}^{m+1}_j}{\Delta t_a},~m=1,2,\cdots,M,
\end{equation}
where $m$ is the sub-iteration step.
When $\Delta \bm{W}^m_j$ converges, the RHS of Eq. (\ref{implicit-2}) becomes a small quantity $\epsilon$, and thus $\bm{W}^{m+1}_j$ is the desired unsteady solution $\bm{W}^{n+1}_j$.
The spatial operator $\mathcal{L}^{m+1}_j$ can be obtained through linearization. Rearranging Eq. (\ref{implicit-2}) yields a linear system with $\Delta \bm{W}^m_s$ as the unknown:
\begin{equation}\label{implicit-3}
\big(\frac{1}{\Delta t}+\frac{1}{\Delta t_a}\big)\triangle \mathbf{W}^m_j-\sum_s\frac{\partial \mathcal{L}^m_j}{\partial \mathbf{W}^m_s} \triangle \mathbf{W}^m_s= \mathcal{L}^{m}_j+\frac{\mathbf{W}^{n}_j-\mathbf{W}^{m}_j}{\Delta t},
\end{equation}
where $s$ represents the cell index for evaluating $\mathcal{L}^m_j$. The derivative $\partial \mathcal{L}^m_j/\partial \mathbf{W}^m_s$ is derived using a simplified method described in the Appendix.
Compared to the original implicit discrete Eq. (\ref{implicit-1}), Eq. (\ref{implicit-3}) introduces an error of $O(\epsilon)/\Delta t_a + O(\epsilon^2)$, where $\epsilon$ is the residual at convergence. Additionally, the evaluation of the Jacobian matrix $\partial\mathcal{L}_j^m/\partial\bm{W}^m_s$ may also introduce errors.
Finally, the implicit scheme based on backward Euler discretization and dual time-stepping method can be obtained as
\begin{equation}\label{implicit-4}
\begin{split}
&\mathbf{A} \triangle \mathbf{W}^m= \mathcal{R}^m,\\
&\mathbf{A}=\big(\frac{1}{\Delta t}+\frac{1}{\Delta t_a}\big)\mathbf{I}-\big(\frac{\mathrm{d}\mathcal{L}}{\mathrm{d}\mathbf{W}}\big)^m, ~\mathcal{R}^m=\mathcal{L}^{m}+\frac{\mathbf{W}^{n}-\mathbf{W}^{m}}{\Delta t}.
\end{split}
\end{equation}
Here the coefficient matrix $\mathbf{A}$ is a block pentagonal matrix assembled from the implicit scheme Eq. (\ref{implicit-3}) over all cells.

Eq. (\ref{implicit-4}) is typically solved using the LU-SGS method. With the LU-SGS method, Eq. (\ref{implicit-4}) is further rewritten in the form of a product of simple lower and upper triangular matrices as follows,
\begin{equation}\label{implicit-LUSGS}
\begin{split}
\big[(\mathbf{L}+\mathbf{D})\mathbf{D}^{-1}(\mathbf{D}+\mathbf{U})\big] \triangle \mathbf{W}^m= \mathcal{R}^m,
\end{split}
\end{equation}
where $\mathbf{L}$, $\mathbf{D}$, and $\mathbf{U}$ represent the strictly lower triangular, diagonal, and strictly upper triangular parts of matrix $\mathbf{A}$.
Exploiting the special structure of the coefficient matrix in Eq. (\ref{implicit-LUSGS}), the equation is decomposed into two systems as follows,
\begin{align*}
\begin{split}
&(\mathbf{L}+\mathbf{D}) \mathbf{y}= \mathcal{R}^m,\\
&\big[\mathbf{D}^{-1}(\mathbf{D}+\mathbf{U})\big] \triangle \mathbf{W}^m= \mathbf{y}.
\end{split}
\end{align*}
The resulting linear systems are then solved explicitly using a forward-backward sweep method.

\section{Numerical examples}

In this section, numerical examples are presented to verify the properties of the fourth-order GENO scheme developed in this study for 3TRD problems, including high-order accuracy performance, bound-preserving properties, and large time-step temporal integration.
The temporal integration employs both the explicit second-order Runge-Kutta (RK) method and the first-order backward Euler method with dual time-stepping method.
The explicit second-order RK method is implemented in a predictor-corrector formulation based on a middle time step. For explicit time integration, the time step $\Delta t$ is restricted to the order of $(\Delta x)^2$, which is commensurate with the spatial error of the fourth-order GENO method. For the first-order implicit method that allows large time step advancement, results obtained by using different time-step sizes are compared.

In the dual-time-stepping method, the inner iteration convergence criterion is set as the reduction of the maximum infinity norm residual of the three conservative variables $W$ by a prescribed order of magnitude. The specific residual order and the pseudo-time step size are provided in the numerical examples.
The computations are performed using 3D grids and the 3D numerical scheme proposed in this study. For 2D problems, three grid layers are used in the third direction.

\subsection{Accuracy test}

\begin{table}
	\begin{center}
		\def\temptablewidth{0.90\textwidth}
		{\rule{\temptablewidth}{0.70pt}}
		\begin{tabular*}{\temptablewidth}{@{\extracolsep{\fill}}c|cc|cc|cc}
			
			mesh size $h$ & Error($T_e$)  & Order & Error($T_i$) & Order & Error($T_r$) & Order \\
			\hline
            1/5  & 1.0140E-005 &       & 1.2736E-005 &      & 1.5951E-005 &             \\
            1/10 & 6.1542E-007 & 4.04  & 7.7084E-007 & 4.05 & 9.6464E-007 & 4.05        \\
            1/20 & 3.7856E-008 & 4.02  & 4.7367E-008 & 4.02 & 5.9244E-008 & 4.02        \\
            1/40 & 2.3438E-009 & 4.01  & 2.9319E-009 & 4.01 & 3.6666E-009 & 4.01        \\
		\end{tabular*}
		{\rule{\temptablewidth}{0.1pt}}
	\end{center}
	\vspace{-1mm} \caption{\label{3d-accuracy-1-1} Accuracy test: $L_1$ errors and convergence orders of electron, ion, and radiation temperatures with the 4th-order GENO scheme. Explicit 2nd-order RK method is adopted with a time step of $\Delta t = 0.1h^2$.}
\end{table}

\begin{table}
	\begin{center}
		\def\temptablewidth{0.90\textwidth}
		{\rule{\temptablewidth}{0.70pt}}
		\begin{tabular*}{\temptablewidth}{@{\extracolsep{\fill}}c|cc|cc|cc}
			
			mesh size $h$ & Error($T_e$) & Order & Error($T_i$) & Order & Error($T_r$) & Order \\
			\hline
            1/5  & 2.2142E-005 &       & 3.2435E-005 &      & 4.7278E-005 &             \\
            1/10 & 1.5308E-006 & 4.04  & 2.2663E-006 & 4.05 & 3.3521E-006 & 4.05        \\
            1/20 & 1.0175E-007 & 4.02  & 1.5106E-007 & 4.02 & 2.2549E-007 & 4.02        \\
            1/40 & 6.4776E-009 & 4.01  & 9.7118E-009 & 4.01 & 1.4563E-008 & 4.01        \\
		\end{tabular*}
		{\rule{\temptablewidth}{0.1pt}}
	\end{center}
	\vspace{-1mm} \caption{\label{3d-accuracy-1-2} Accuracy test: $L_{\infty}$ errors and convergence orders of electron, ion, and radiation temperatures with the 4th-order GENO scheme. Explicit 2nd-order RK method is adopted with a time step of $\Delta t = 0.1h^2$.}
\end{table}

An unsteady linear 3TRD problem \cite{yu2019finite} is solved on a 3-D computational domain to test the accuracy of the numerical scheme developed in this study.
A degenerate linear problem is employed for the present accuracy test. The parameters in the 3TRD Eq. (\ref{3T-pde}) are set as $c_{\alpha}=1$ $\alpha \in \{e,i,r\}$, $\omega_{i}=1$, and $\omega_{r}=1$.
The analytical solution for this problem is given by
\begin{align*}
T_e=e^t(x^2+1)(y^2+1),~T_i=e^t(2x^2+1)(y^2+1),~T_r=e^t(2x^2+1)(2y^2+1).
\end{align*}
For simplicity in handling boundary conditions, values at ghost cells are prescribed by the analytical solution.
To satisfy the given solution, an additional source term is required as
\begin{align*}
\begin{split}
S_e^{*}&=-e^t(3x^2y^2+3x^2+2y^2+3),\\
S_i^{*}&=e^t(3x^2y^2-x^2-3y^2-5),\\
S_r^{*}&=e^t(7x^2y^2-5x^2-5y^2-7).
\end{split}
\end{align*}
The computational domain is set to $[0,1]^2 \times [0,3h]$, where $h$ denotes the grid spacing. For the 2D problem, the grid is maintained at three cells in the third direction.
The simulation is run until $t = 1$ with a time step of $\Delta t = 0.1(\Delta x)^2$.

\begin{figure}[!htb]
\centering
\includegraphics[width=0.45\textwidth]{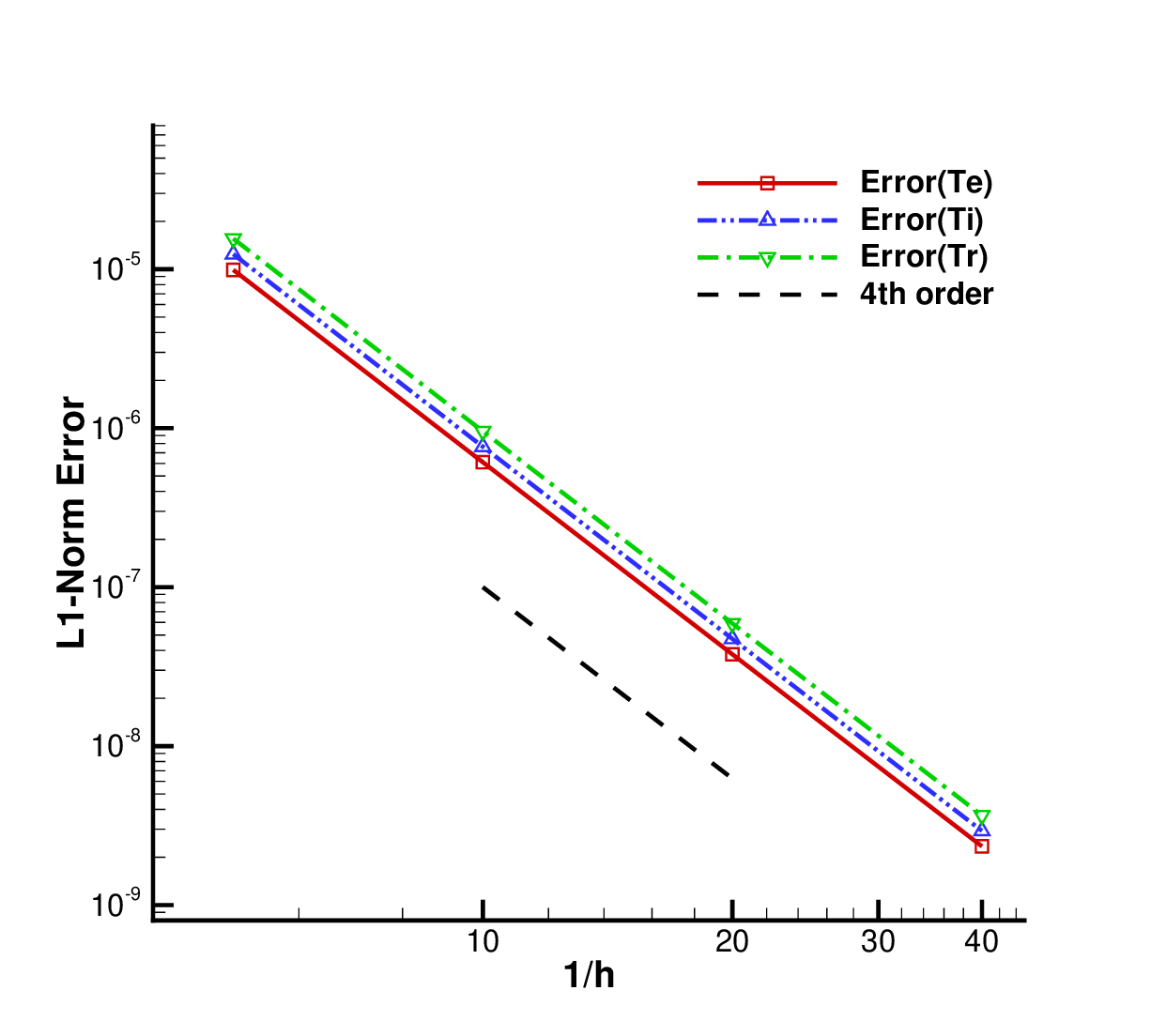}
\includegraphics[width=0.45\textwidth]{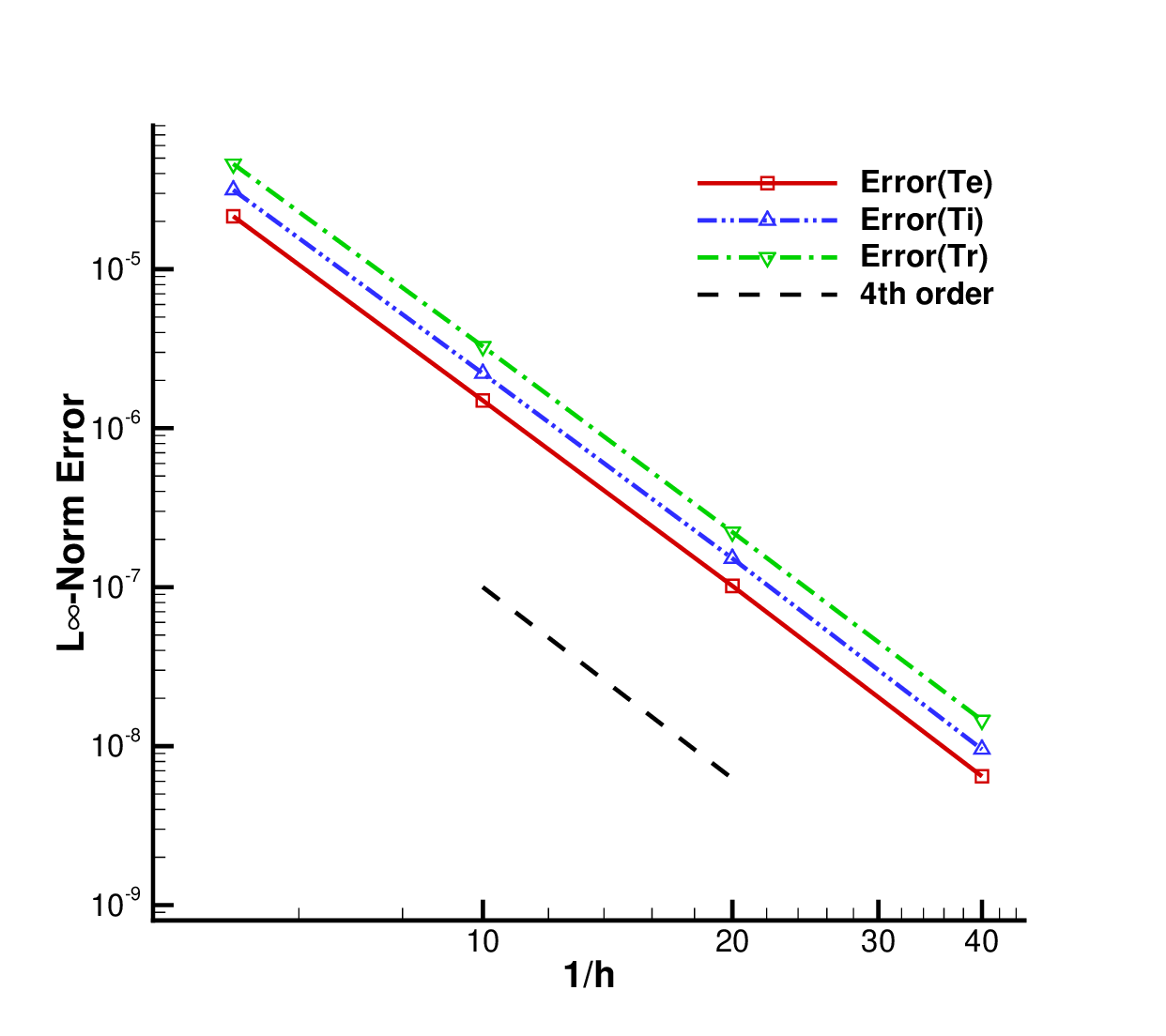}
\caption{\label{3d-accuracy-conver-line} Accuracy test: Convergence of the $L_1$ and $L_\infty$ error norms for the electron ($T_e$), ion ($T_i$), and radiation ($T_r$) temperatures under mesh refinement. }
\end{figure}

Tables \ref{3d-accuracy-1-1} and \ref{3d-accuracy-1-2} present the $L_1$ and $L_\infty$ error norms, respectively, along with the convergence orders for the three temperatures as the grid is refined.
The results demonstrate that the present GENO method achieves the theoretical fourth-order accuracy. Figure \ref{3d-accuracy-conver-line} illustrates the error reduction with grid refinement.

\subsection{2D model problem}

This test case is employed to verify the bound-preserving property and high-accuracy performance of the fourth-order GENO scheme in the presence of temperature discontinuities and local large gradients, as well as to assess the performance of the implicit large-time-step integration method in handling strong stiffness. This test case has also been studied in \cite{yu2019finite,peng2020positivity}.
The computational domain $\Omega = [0,300]^2 \times [0,3h]$ is partitioned into two sub-regions:
\begin{equation*}
\begin{split}
&\Omega_A=\{(x,y,z)\in \Omega: x\in[0,300],y\in[0,250]\},\\
&\Omega_B=\{(x,y,z)\in \Omega: x\in[0,300],y\in[250,300]\},
\end{split}
\end{equation*}
where $h$ is the grid spacing and the grid is maintained at three cells in the third direction. $\Omega_A$ and $\Omega_B$ consist of different materials, with discontinuous material properties across the interface.
Material properties are prescribed as follows:
\begin{equation*}
\begin{split}
&\Omega_A: k_e=k_i=10,~k_r=100,~c_e=c_i=c_r=0.05,~ w_i=10,w_r=100,\\
&\Omega_B: k_e=k_i=k_r=10,~c_e=c_i=c_r=1,~w_i=10,w_r=100.
\end{split}
\end{equation*}
The initial temperature is initialized uniformly to $T_{\alpha}=3\times 10^{-4}$ for all species $\alpha \in \{e,i,r\}$. Regarding boundary conditions, Neumann condition $\partial T_\alpha/\partial n=0$ is applied on all boundaries, with the exception of the left boundary ($x=0$), where a fixed radiation temperature $T_r=100$ is imposed.
A uniform grid with spacing $\Delta x=\Delta y=3$ is employed.
At this mesh spacing, the stiffness of the system restricts the time-step size of the explicit second-order RK method to $3\times 10^{-4}$.
The simulation is terminated at $t=5$.

\subsubsection{Bound-preserving property}

\begin{figure}[!htb]
\centering
\includegraphics[width=0.49\textwidth]{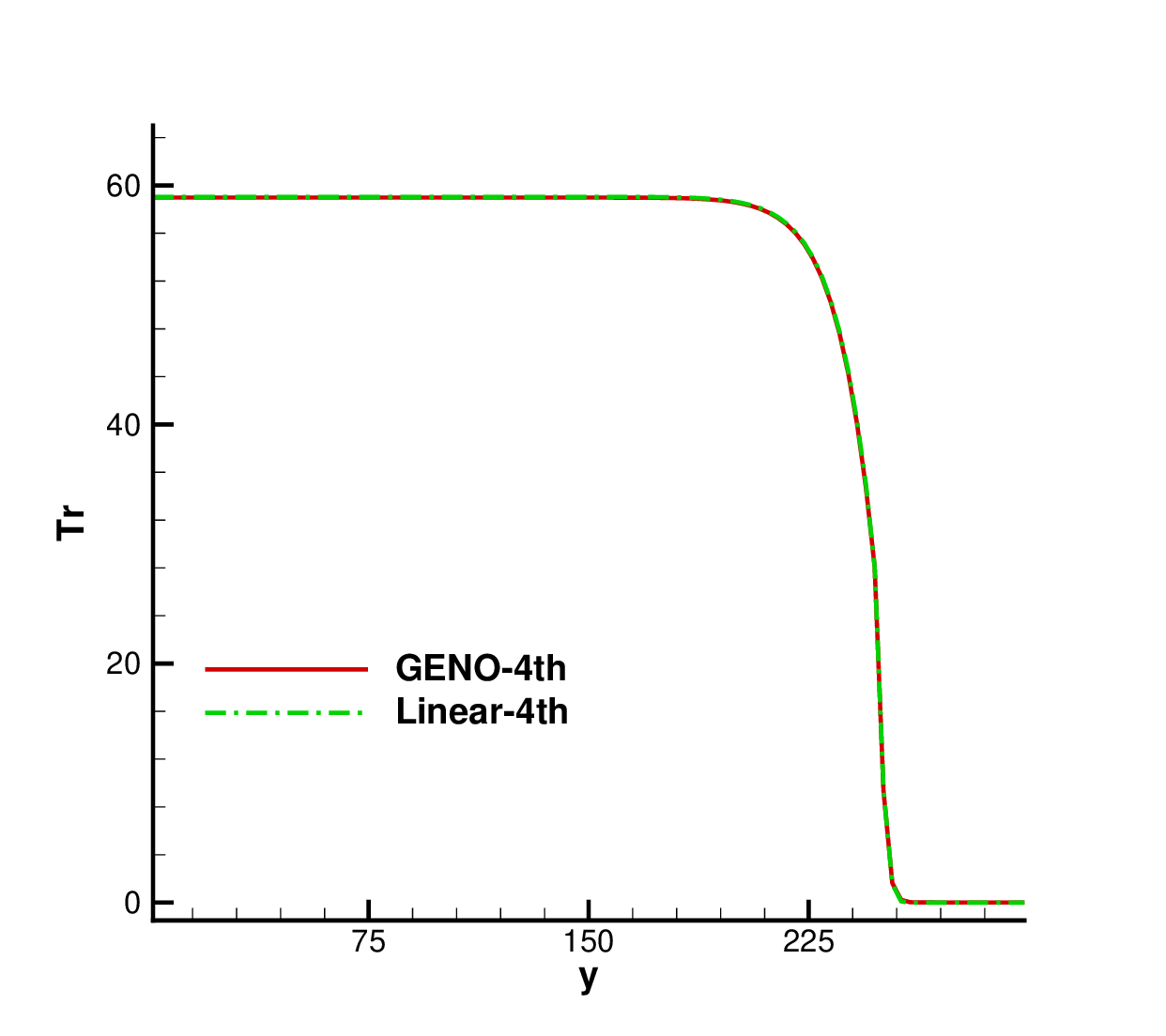}
\includegraphics[width=0.49\textwidth]{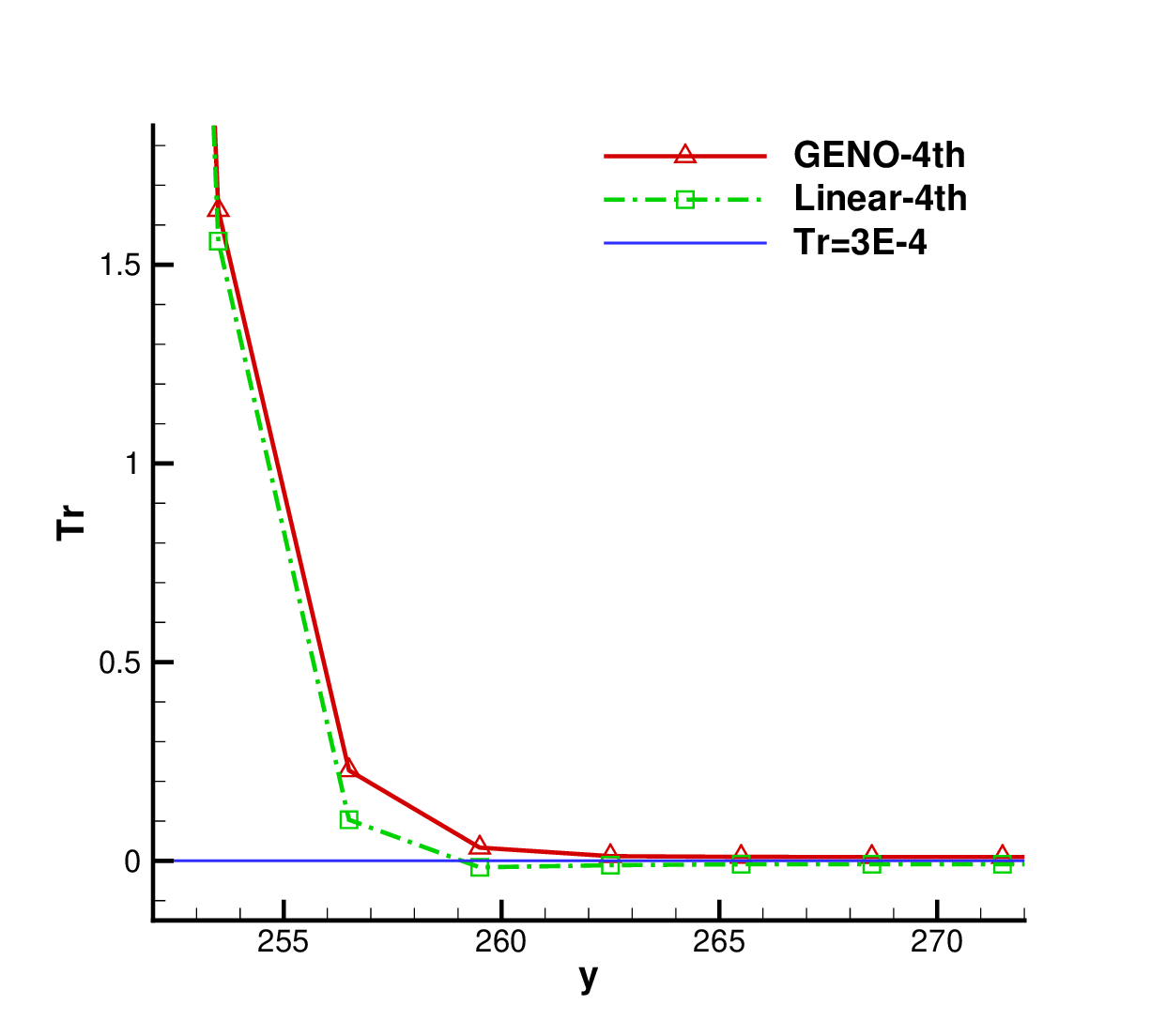}
\caption{\label{2d-mono-test-1} 2D model problem for verifying the bound-preserving property: Radiation temperature distributions along $x=0.15$ at $t=0.5$ computed by the 4th-order GENO and 4th-order linear schemes, with a close-up view across the material interface shown on the right.}
\end{figure}

We first verify the bound-preserving property of the fourth-order GENO scheme. Due to the presence of a very large initial temperature jump at the left boundary, as well as the strong discontinuities in material parameters on the two sides of $y=250$, extremely large temperature gradients are generated, which pose a significant challenge to high-order numerical schemes.
Table \ref{2d-mono-test-table} lists the temperature bounds computed at an early stage ($t=0.5$) using the fourth-order scheme with and without GENO reconstruction; the time discretization is performed using the explicit second-order RK method.
The results demonstrate that the GENO scheme strictly preserves physical bounds, whereas the standard fourth-order linear scheme yields nonphysical negative temperatures.
Additionally, Figure \ref{2d-mono-test-1} displays the temperature profile at $x=0.15$ across the material interface. The steep gradient at the interface causes the linear scheme to oscillate and produce negative values, a problem effectively resolved by the GENO reconstruction.

\begin{table}[htbp]
    \centering
    \begin{tabular}{c|c|c|c|c|c|c}
        \hline
        Reconstruction       & $T_{e,min}$ & $T_{e,max}$ & $T_{i,min}$ & $T_{i,max}$ & $T_{r,min}$ & $T_{r,max}$ \\ \hline
        GENO-4th             & $3.0\times 10^{-4}$ & 94.9985 & $3.0\times 10^{-4}$ & 94.1166 & $3.0\times 10^{-4}$ & 95.1827 \\ \hline
        Linear-4th           & $-0.0303438$ & 94.9969 & $-0.0265335$ & 94.0719 & $-0.0305249$ & 95.1921 \\ \hline
    \end{tabular}
    \caption{2D model problem for verifying the bound-preserving property: Temperature bounds at $t=0.5$ computed by the 4th-order schemes using GENO and linear reconstructions. The GENO scheme successfully preserves the physical bounds.}
    \label{2d-mono-test-table}
\end{table}

\begin{figure}[!htb]
\centering
\includegraphics[width=0.49\textwidth]{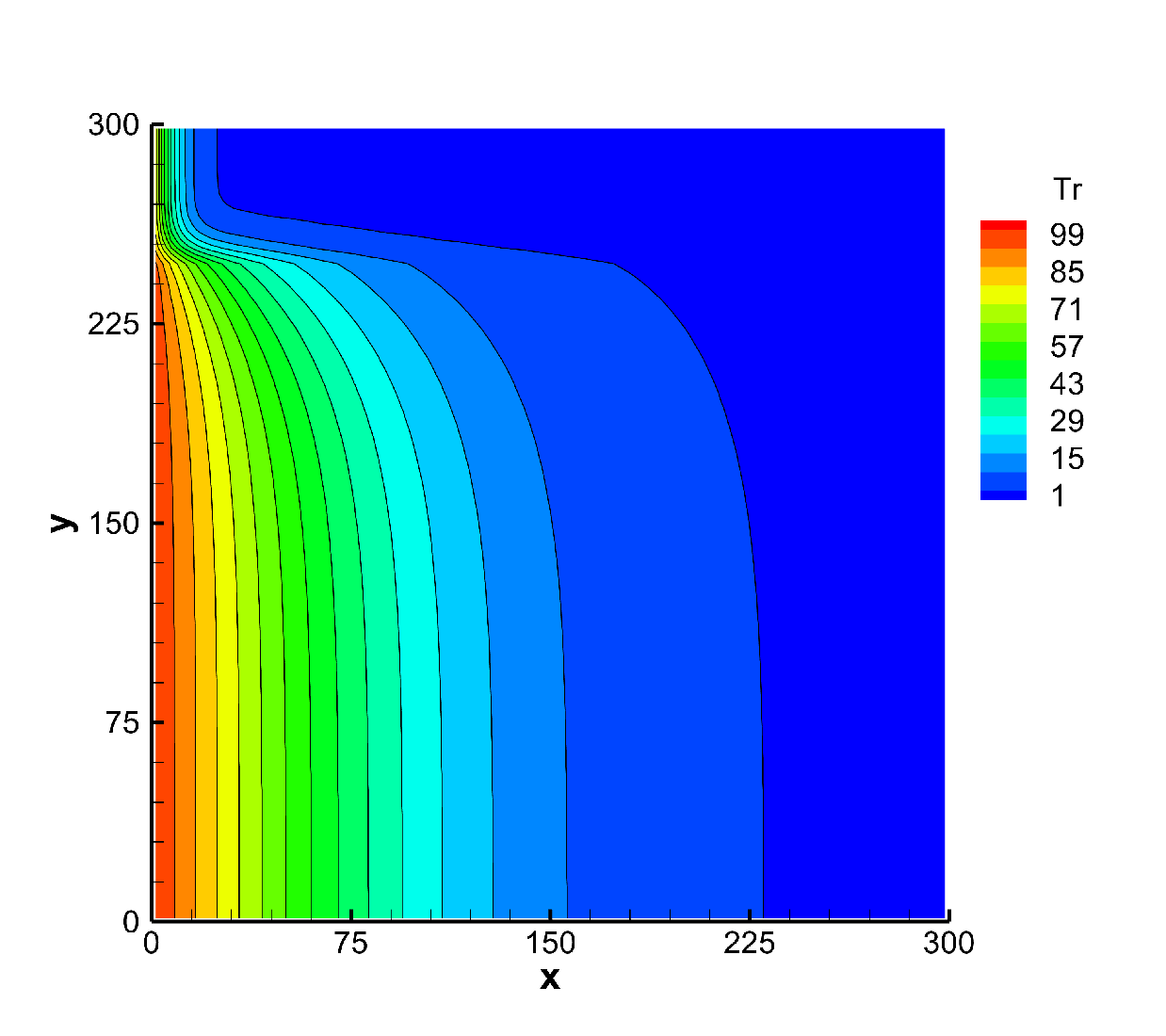}
\includegraphics[width=0.49\textwidth]{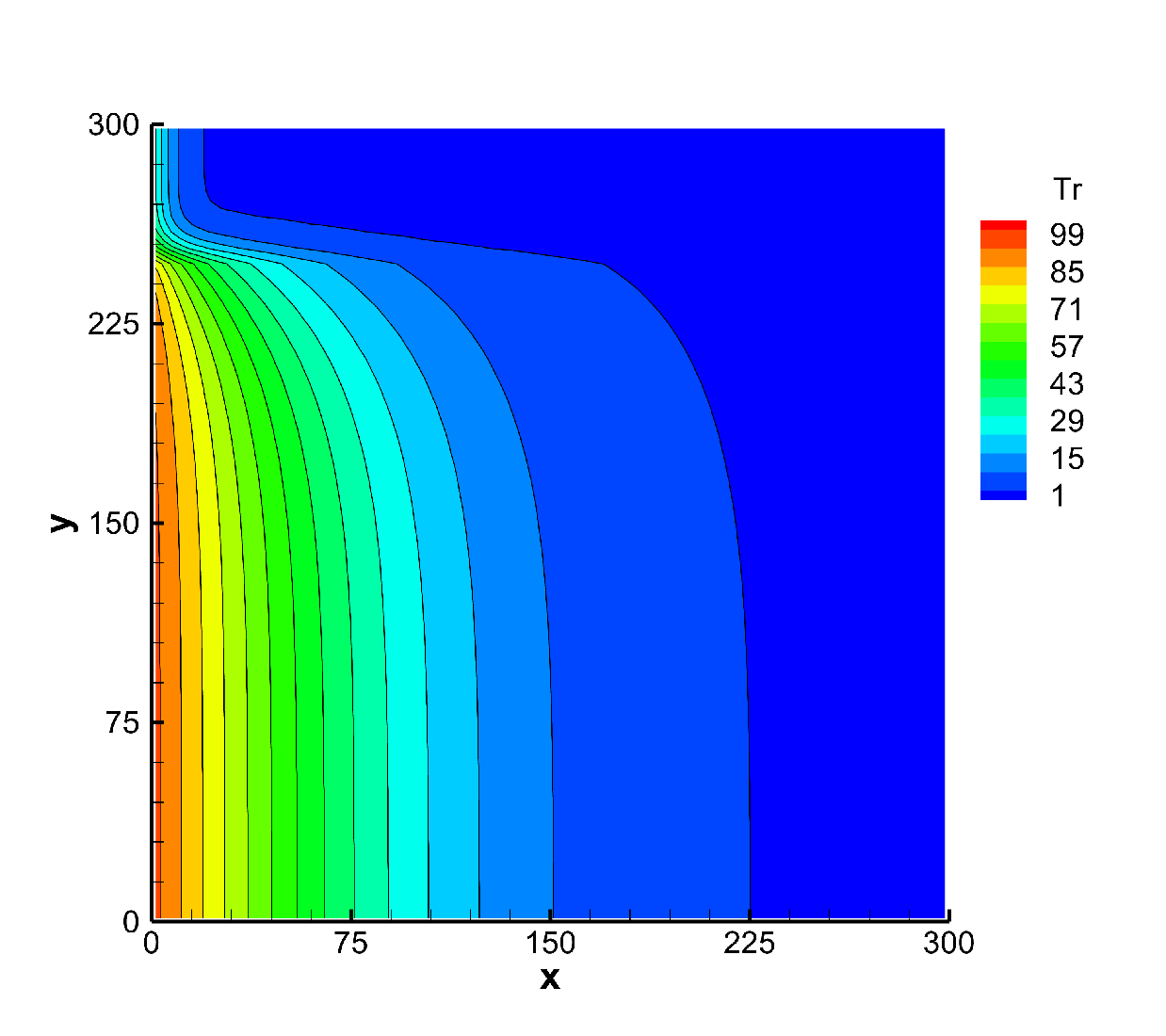}
\caption{\label{2d-mono-test-2} 2D model problem for evaluating high accuracy performance: Contours of $T_r$ at $t=5$ computed by the 4th-order GENO scheme (left) and the linear 2nd-order central scheme (right).}
\end{figure}

\begin{figure}[!htb]
\centering
\includegraphics[width=0.49\textwidth]{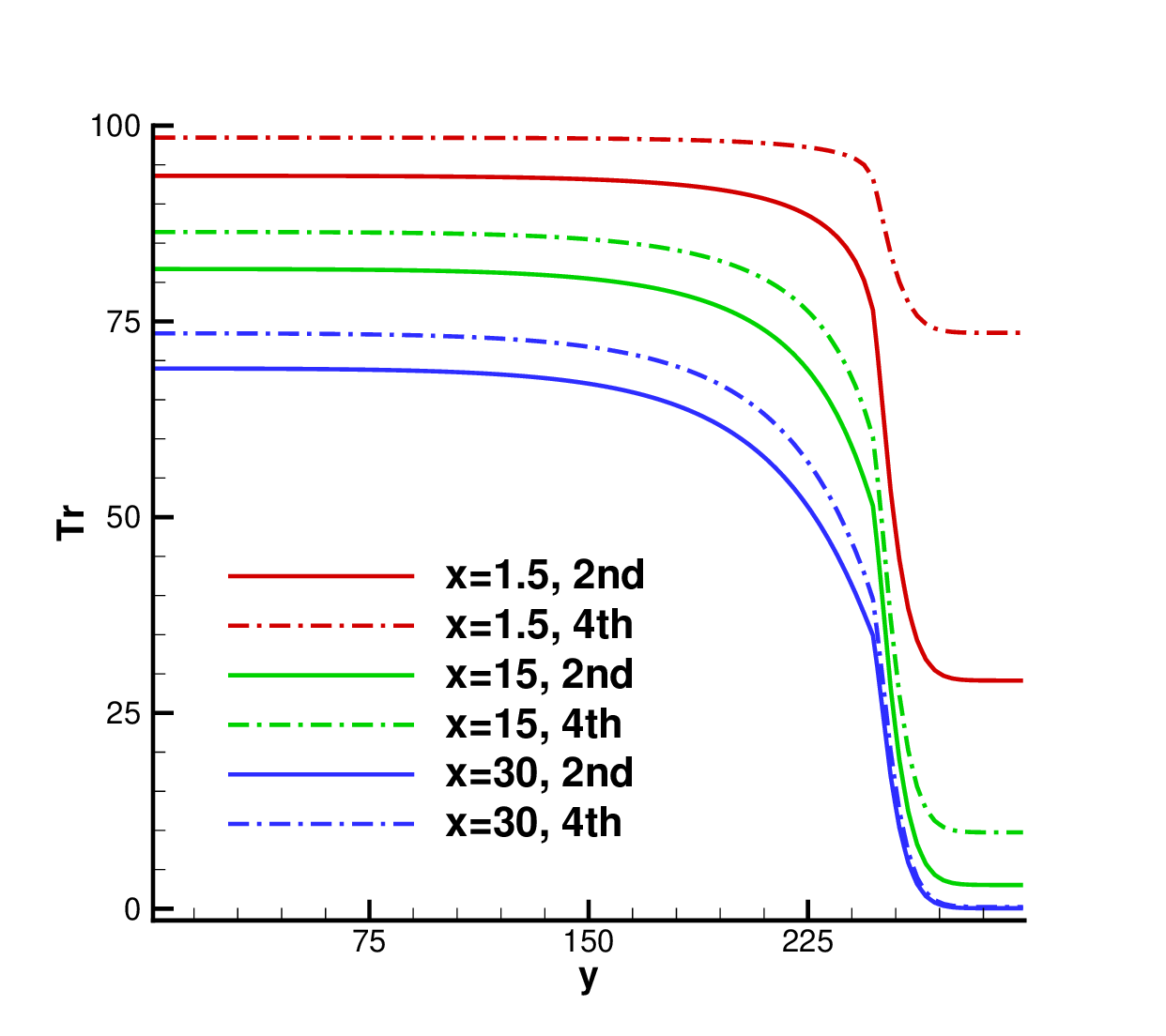}
\includegraphics[width=0.49\textwidth]{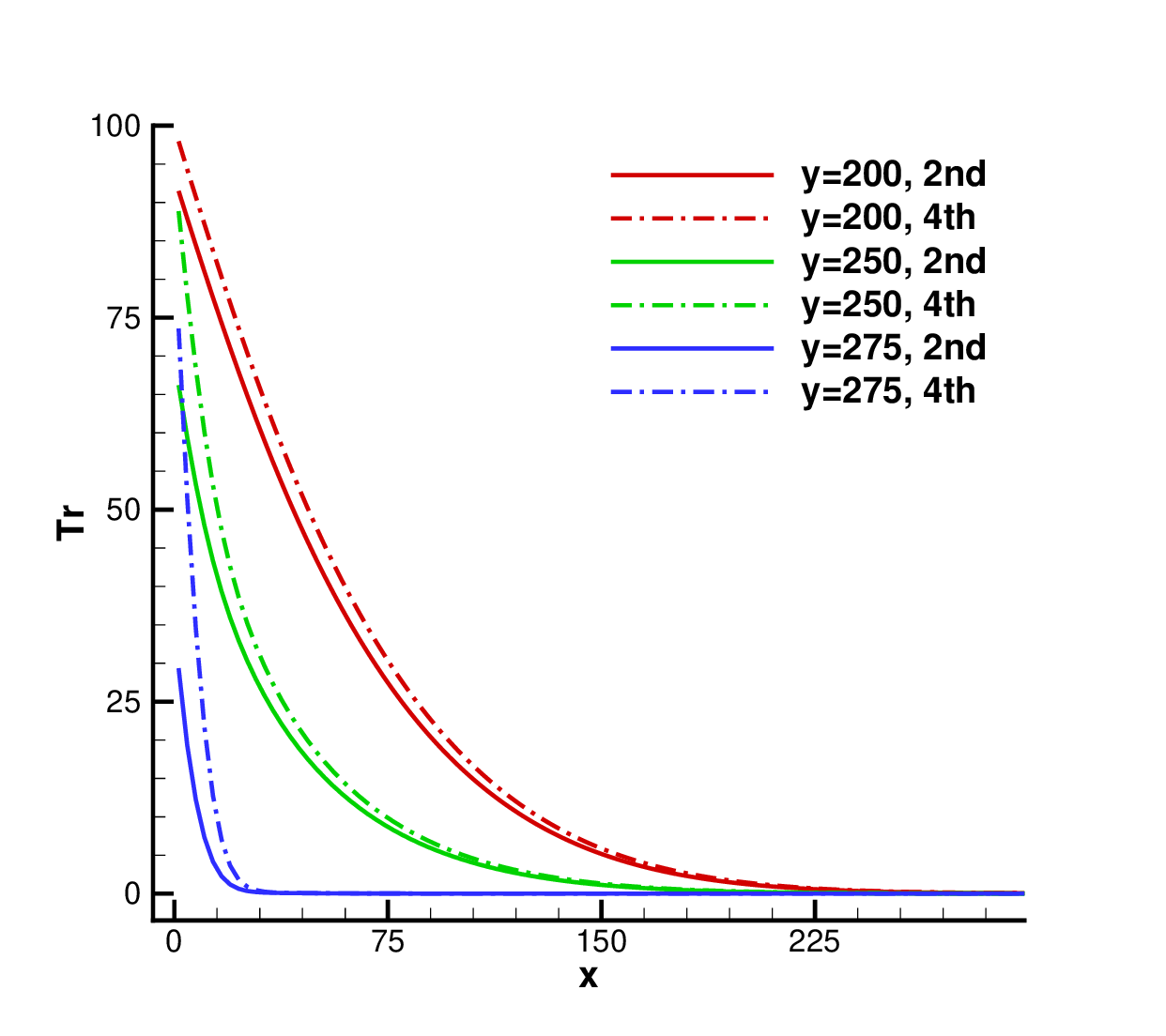}
\caption{\label{2d-mono-test-3} 2D model problem for evaluating high accuracy performance: Quantitative comparison of \(T_r\) at \(t=5\) computed by the 4th-order GENO scheme and the linear 2nd-order central scheme along various lines in the \(x\)-direction (left) and \(y\)-direction (right).}
\end{figure}

\subsubsection{High accuracy performance}

Furthermore, this test case is also used to compare the accuracy performance of fourth-order and second-order schemes for the 3TRD problem, with time discretization uniformly employing the second-order explicit RK method. Figure \ref{2d-mono-test-2} presents the computational results using fourth-order GENO reconstruction and second-order linear reconstruction, where the left and right panels show the radiation temperature contours for the fourth-order and second-order schemes, respectively.
Near the left boundary of the computational domain, the fourth-order scheme yields higher temperatures, and overall, thermal energy is transported over greater distances. Figure \ref{2d-mono-test-3} provides a quantitative comparison of temperature distributions at different locations, where the fourth-order scheme produces higher temperature distributions.
This behavior is attributed to the presence of a steep temperature gradient at the left boundary of the computational domain. The fourth-order spatial reconstruction resolves this gradient more accurately, thereby ensuring a more precise, and consequently higher, heat flux transport. This highlights the accuracy advantage of the high-order scheme for this problem.

\begin{figure}[!htb]
\centering
\includegraphics[width=0.49\textwidth]{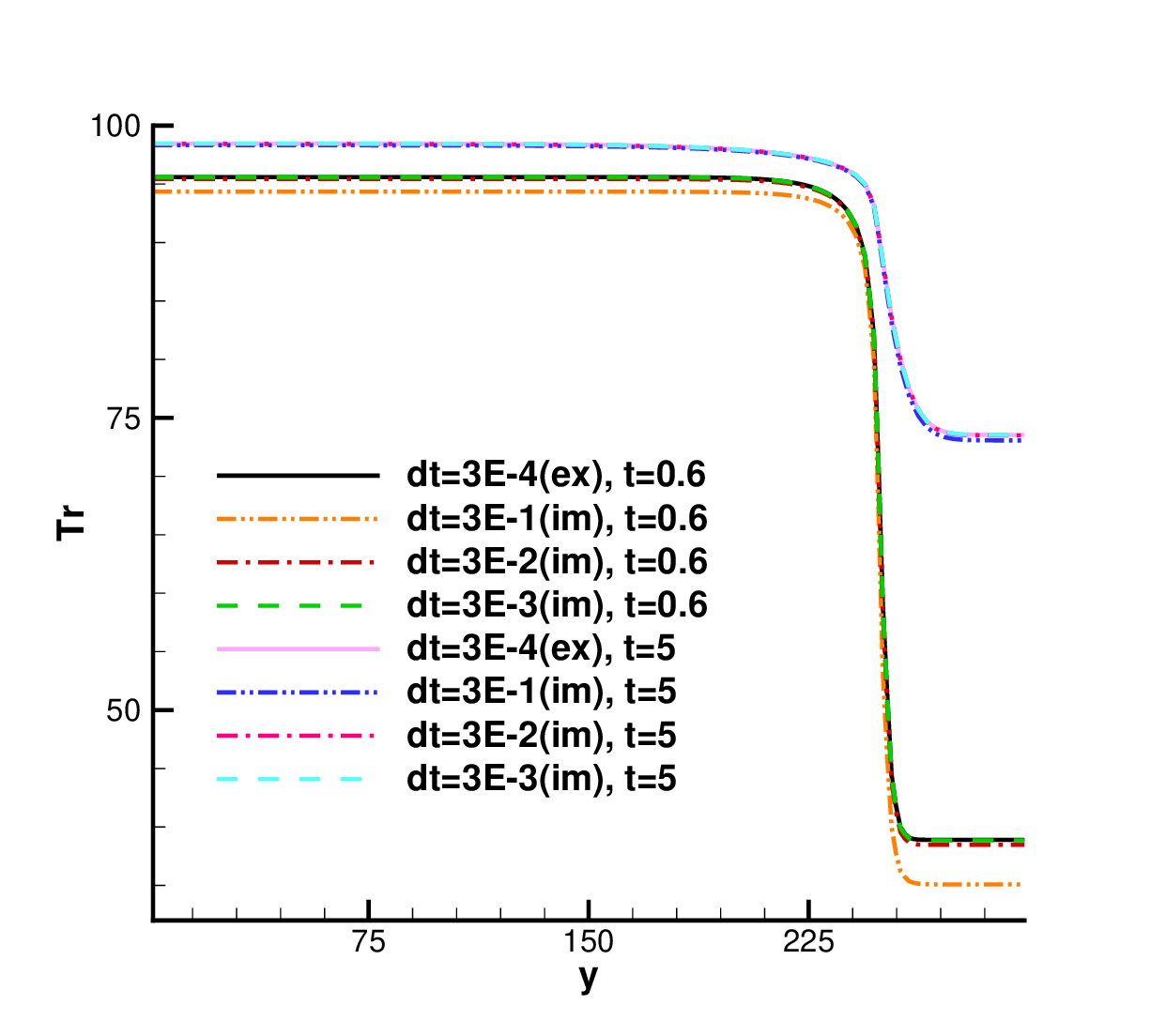}
\includegraphics[width=0.49\textwidth]{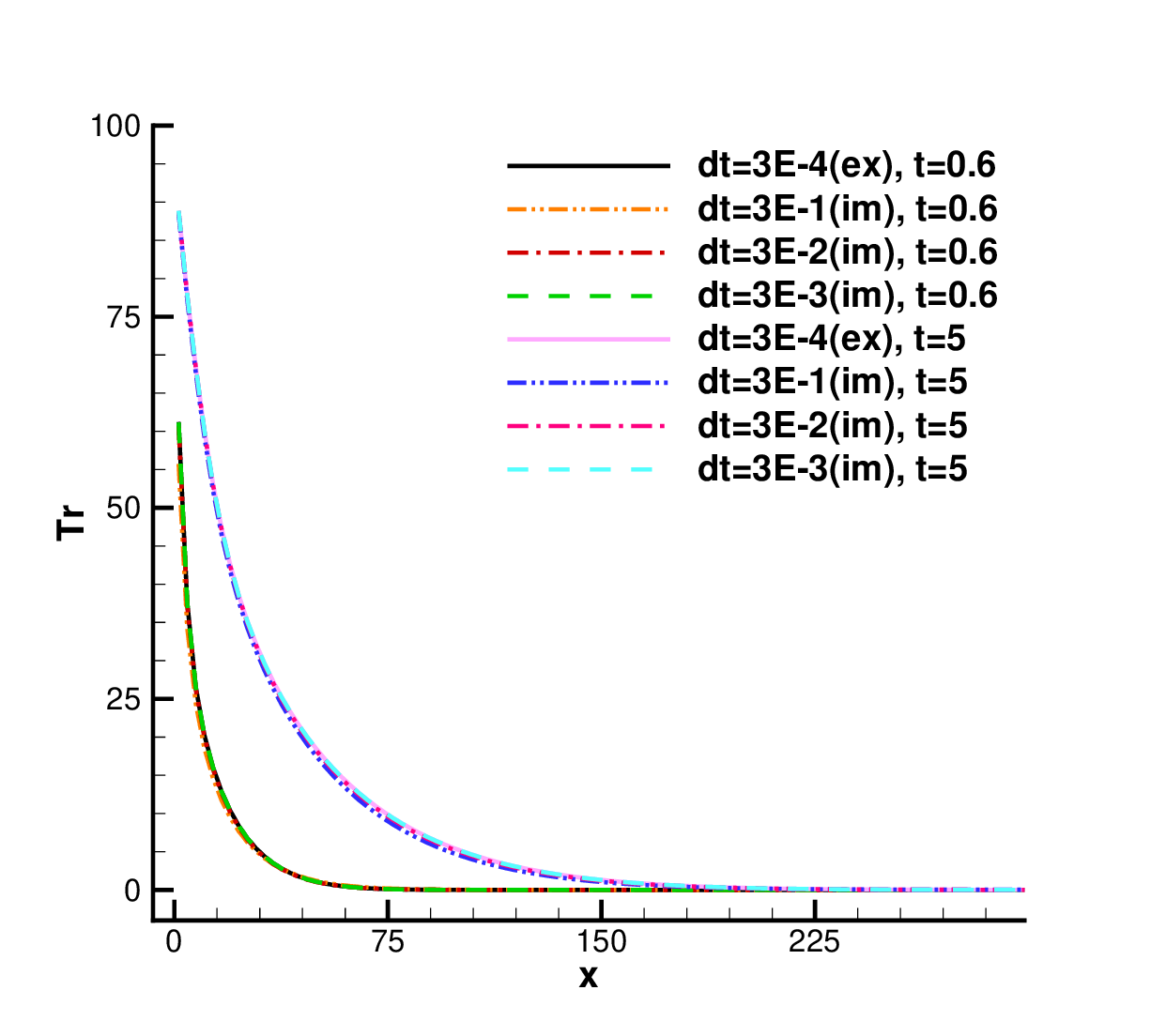}
\caption{\label{2d-mono-test-4} 2D model problem for evaluating the implicit large time-step integration: Quantitative comparison of temperature distributions along various lines in the \(x\)-direction (left) and \(y\)-direction (right) computed using different time-step sizes.}
\end{figure}

\subsubsection{Implicit large time-step integration}

Finally, we verify the implicit large time-step integration method. Figure \ref{2d-mono-test-4} shows quantitative comparisons of temperature distributions at different times and positions for time steps of 10$\Delta t$, 100$\Delta t$, and 1000$\Delta t$, where $\Delta t=3\times 10^{-4}$ represents the maximum time step permitted by the explicit RK method.
The inner iteration convergence criterion is set such that the largest $L_\infty$ norm of the three temperature residuals decreases by $4$ orders of magnitude. Additionally, the pseudo time step is set to $\Delta_a=1000\Delta t$.
Since the maximum temperature gradients occur at $x=1.5$ and $y=250$, these locations provide a more stringent accuracy test. The results demonstrate that time steps of 10$\Delta t$ and 100$\Delta t$ yield solutions nearly identical to the explicit small time-step solution, with slight deviations only appearing at 1000$\Delta t$. This confirms that the proposed implicit time discretization enables large time steps while preserving computational accuracy.

\subsection{Simplified ICF model problem}
We consider the typical model for laser-driven implosion in ICF reported in \cite{zeyao2004parallel} as a benchmark test.
The original benchmark problem is posed on a two-dimensional semicircular domain. To facilitate the
implementation and verification of the proposed high-order GENO scheme on Cartesian grids, we instead consider a
simplified 3D cuboidal computational domain \(\Omega\), defined by $\Omega=[-115,115]\times[0,115]\times[0,115]$.
The boundary \(\partial \Omega\) is divided into two parts. The planes \(y=0\) and \(z=0\) are symmetry boundaries,
denoted by \(\Gamma_1\), on which Neumann boundary conditions are imposed for all three species as $\partial T_{\alpha}/\partial n=0$ ($\alpha \in \{e,i,r\}$).
The remaining boundary faces constitute the laser-irradiated isothermal boundary, denoted by \(\Gamma_2\): an
isothermal condition is prescribed for radiation temperature as $T_r=2.0$, whereas adiabatic
boundary conditions are applied to $T_{e}$ and $T_{i}$. The initial temperature is initialized uniformly to $T_{\alpha}=3.0 \times10^{-4}$ for three species.
Analogously, the domain is partitioned into three subregions (inner, middle, and outer layers), defined as
\begin{align*}
\begin{split}
&\Omega_1=\left\{(x,y,z)\in \Omega:\ |x|\in[0,85],\ y\in[0,85],\ z\in[0,85]\right\},\\
&\Omega_3=\left\{(x,y,z)\in \Omega:\ |x|\in[95,115]\ \text{or}\ y\in[95,115]\ \text{or}\ z\in[95,115]\right\},\\
&\Omega_2 = \Omega\setminus\left(\Omega_1\cup \Omega_3\right).
\end{split}
\end{align*}
The three subregions are filled with deuterium gas (\(DT\)), glass (\(SiO_2\)), and plastic foam (\(CH\)), respectively. All parameters are set as follows:
\begin{align*}
&c_{v\alpha}=
\begin{cases}
1.5\,\Gamma_e, & \alpha=e,\\
1.5\,\Gamma_i, & \alpha=i,\\
\Gamma_r,      & \alpha=r,
\end{cases}
\qquad
k_{\alpha}=
\begin{cases}
A_eT_e^{5/2},     & \alpha=e,\\
A_iT_i^{5/2},     & \alpha=i,\\
A_rT_r^{\beta+3}, & \alpha=r,
\end{cases}
\qquad
\omega_{\alpha}=
\begin{cases}
\rho^2 A_{ei}T_e^{-2/3},     & \alpha=i,\\
\rho^2 A_{er}T_e^{-1/2},     & \alpha=r,
\end{cases} \\
&\rho=
\begin{cases}
0.09, & \text{in } \Omega_1,\\
2.50, & \text{in } \Omega_2,\\
1.10, & \text{in } \Omega_3,
\end{cases}
\qquad
\Gamma_e=
\begin{cases}
35, & \text{in } \Omega_1,\\
40, & \text{in } \Omega_2,\\
45, & \text{in } \Omega_3,
\end{cases}
\qquad
\Gamma_i=
\begin{cases}
35, & \text{in } \Omega_1,\\
40, & \text{in } \Omega_2,\\
70, & \text{in } \Omega_3,
\end{cases}
\qquad
\Gamma_r=0.007568,\\
&A_e=
\begin{cases}
200, & \text{in } \Omega_1,\\
60,  & \text{in } \Omega_2,\\
81,  & \text{in } \Omega_3,
\end{cases}
\qquad
A_i=
\begin{cases}
5,              & \text{in } \Omega_1,\\
1.7\times10^{-4}, & \text{in } \Omega_2,\\
2.0\times10^{-2}, & \text{in } \Omega_3,
\end{cases}
\qquad
A_r=
\begin{cases}
1.8\times10^{7}/\rho,        & \text{in } \Omega_1,\\
9.0\times10^{2}/\rho^{3/2},  & \text{in } \Omega_2,\\
2.1\times10^{3}/\rho^{2},    & \text{in } \Omega_3,
\end{cases}\\
&\beta=
\begin{cases}
1.0, & \text{in } \Omega_1,\\
2.4, & \text{in } \Omega_2,\\
3.0, & \text{in } \Omega_3,
\end{cases}
\qquad
A_{ei}=
\begin{cases}
2000, & \text{in } \Omega_1,\\
4000, & \text{in } \Omega_2,\\
7000, & \text{in } \Omega_3,
\end{cases}
\qquad
A_{er}=
\begin{cases}
10,  & \text{in } \Omega_1,\\
140, & \text{in } \Omega_2,\\
79,  & \text{in } \Omega_3.
\end{cases}
\end{align*}

The computations are performed on a uniform grid with mesh spacing $h=5$. The time-step size is set to $\Delta t=6\times 10^{-6}$.
For this problem, explicit time-marching methods are computationally prohibitive owing to the extremely small time-step size imposed by the severe stiffness;
therefore, only the implicit method is employed. The pseudo-time step for the dual time-stepping method is taken as $\Delta\tau=10\Delta t$, and the inner iterations are considered converged
when the residual has been reduced by $3$ orders of magnitude.

\begin{figure}[!htb]
\centering
\includegraphics[width=0.475\textwidth]{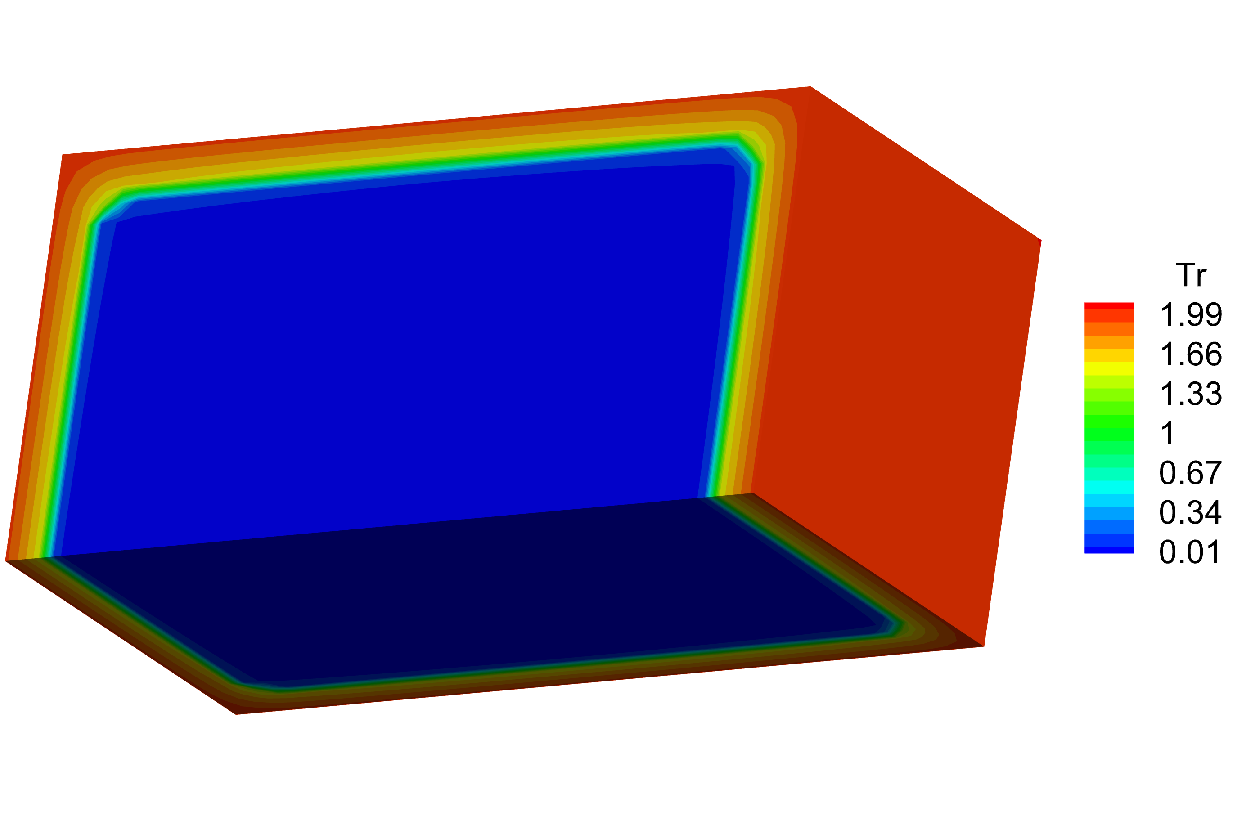}
\includegraphics[width=0.475\textwidth]{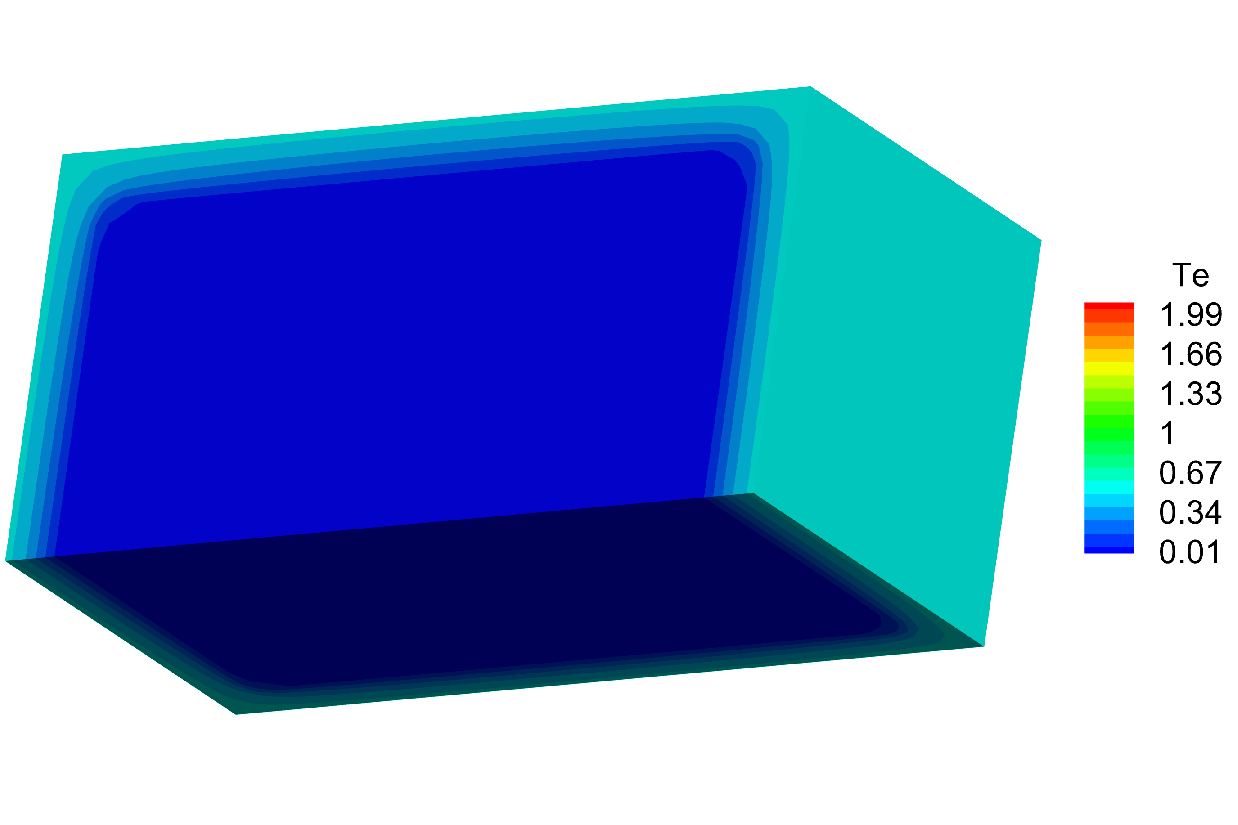}\\
\includegraphics[width=0.475\textwidth]{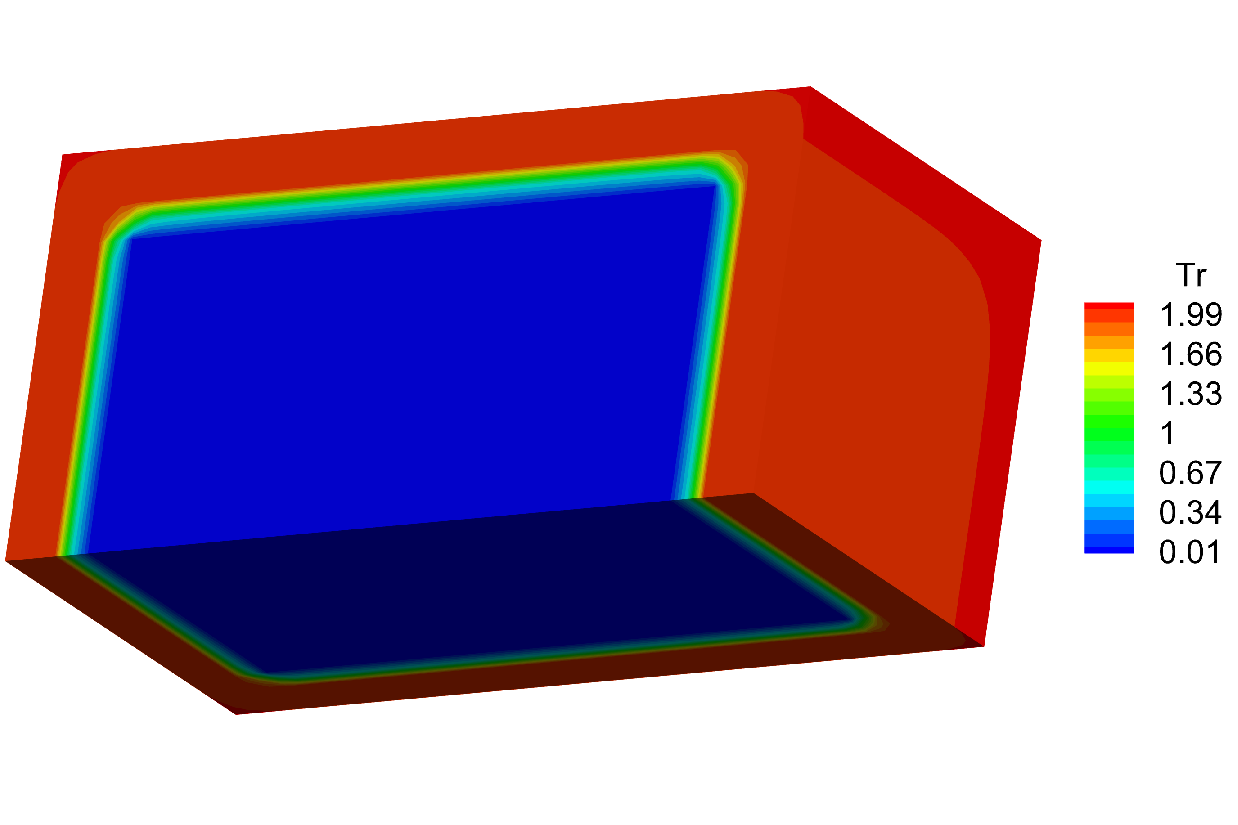}
\includegraphics[width=0.475\textwidth]{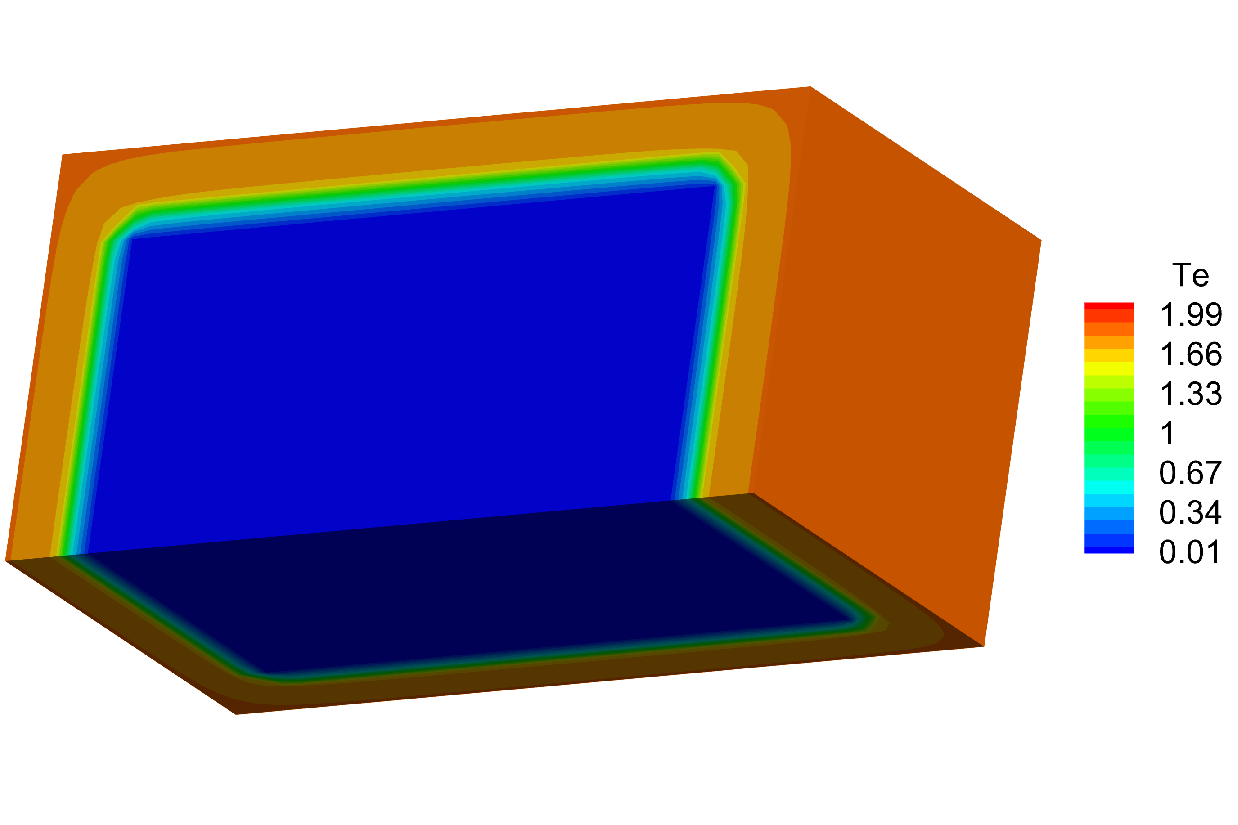}\\
\includegraphics[width=0.475\textwidth]{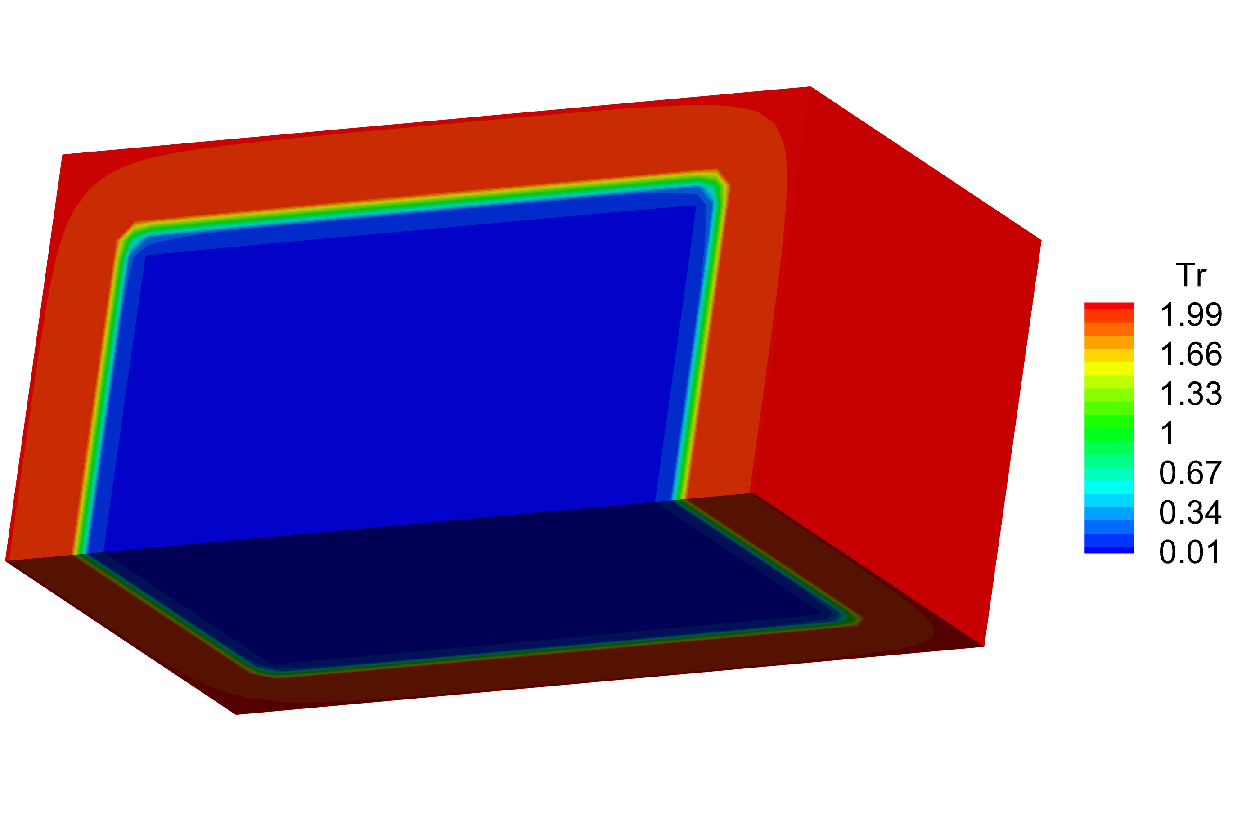}
\includegraphics[width=0.475\textwidth]{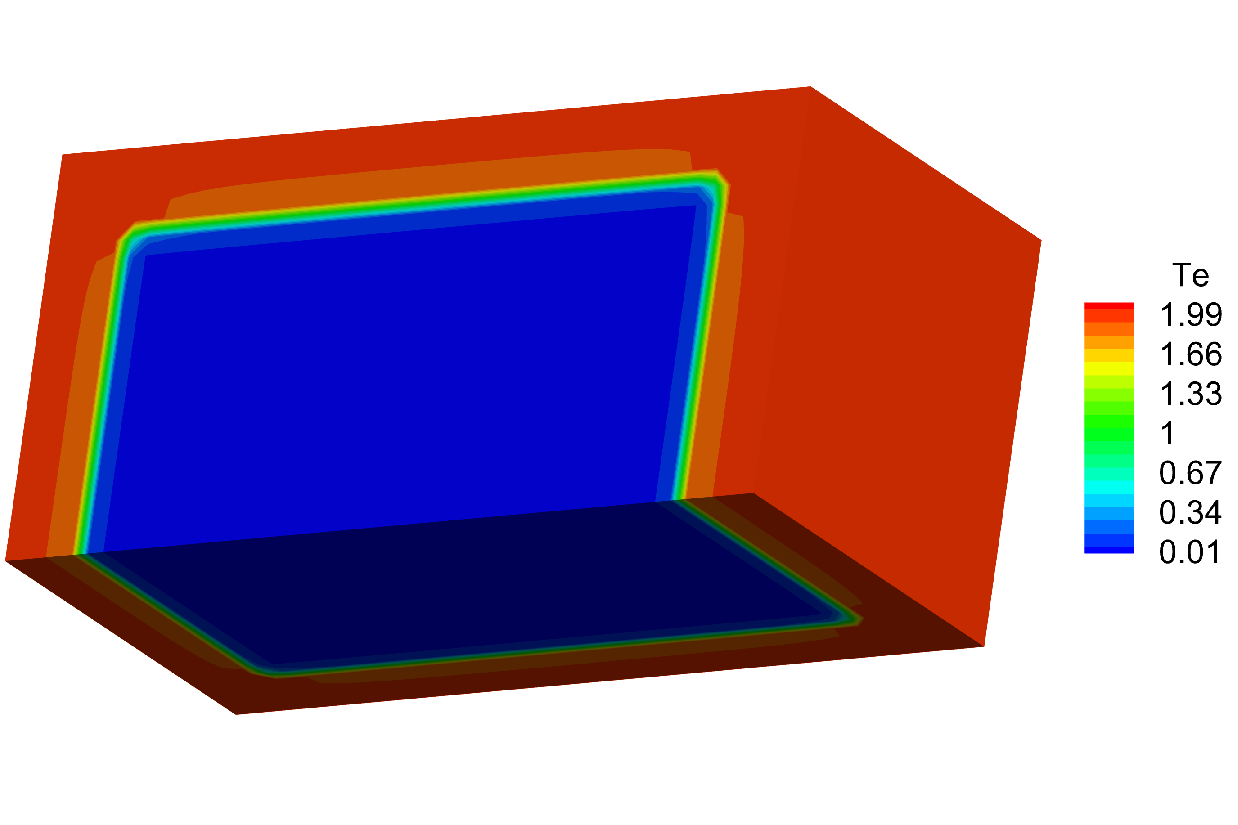}
\caption{\label{3d-ICF-model-1} 3D ICF model problem: Contours of \(T_r\) (left) and \(T_e\) (right) at \(t=0.3\) (top), \(t=5.0\) (middle), and \(t=10\) (bottom), computed using the 4th-order GENO scheme with implicit time integration.}
\end{figure}

Figure \ref{3d-ICF-model-1} presents the contour plots of $T_r$ and $T_e$ obtained using the fourth-order GENO scheme.
Since the electron temperature increases primarily through energy exchange with radiation, its rise is relatively slow, approaching the wall temperature of $2$ only after an extended period of evolution.
Table \ref{3d-ICF-model-table} presents the lower and upper bounds of the three temperatures at simulation times $t=0.3$ and $t=5.0$.
The results demonstrate that all three temperatures strictly satisfy the bound-preserving property.

\begin{figure}[!htb]
\centering
\includegraphics[width=0.49\textwidth]{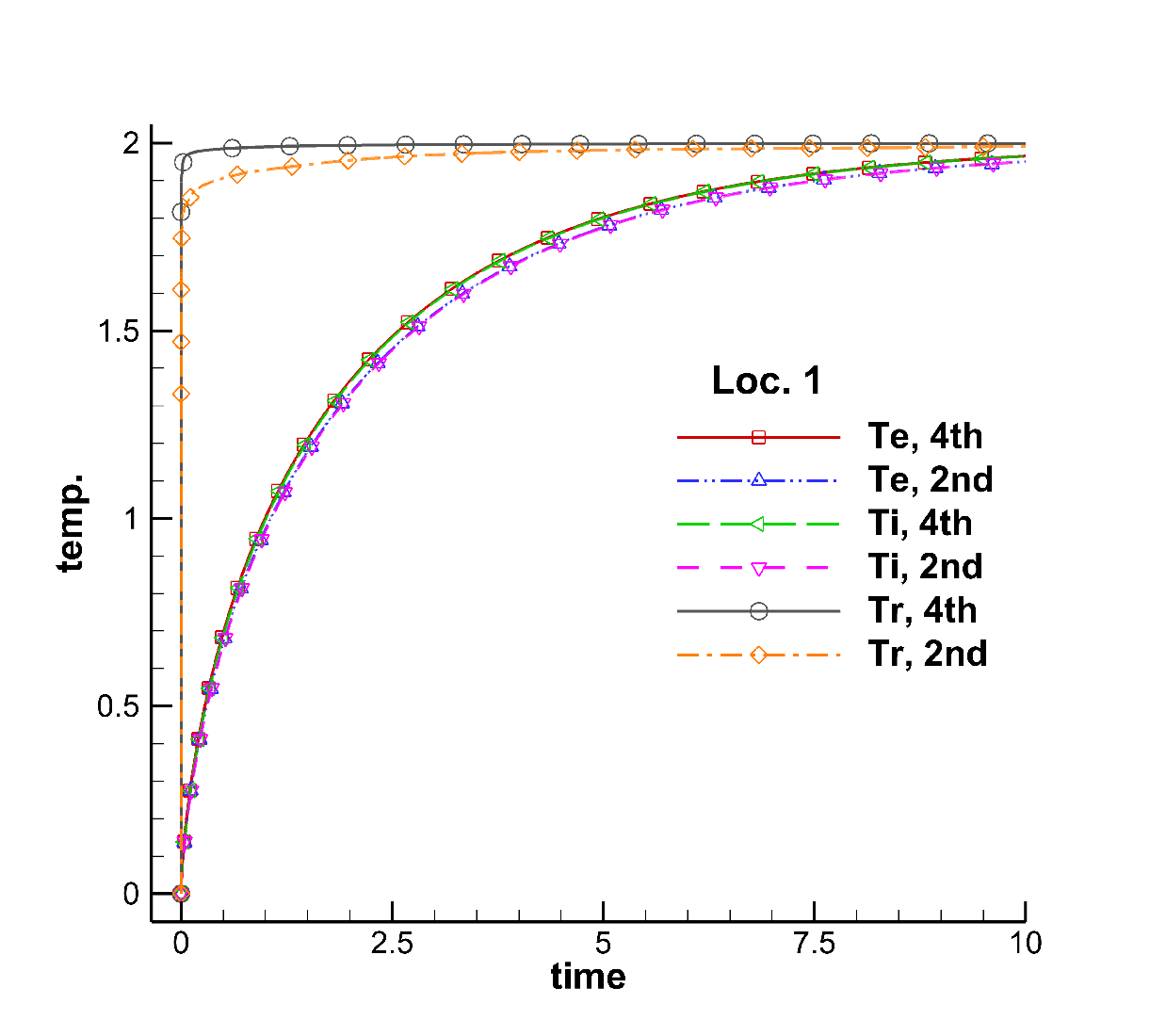}
\includegraphics[width=0.49\textwidth]{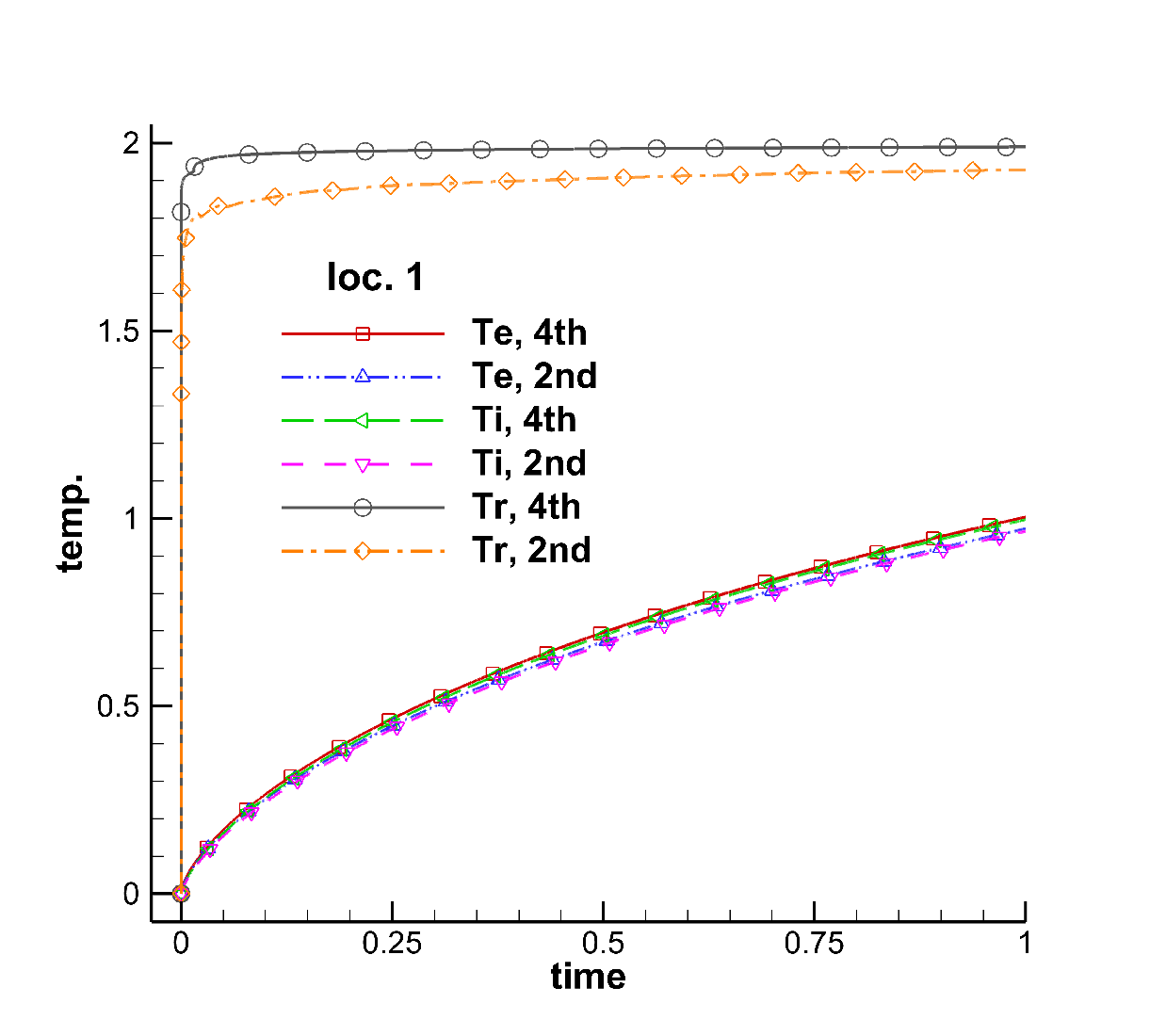}\\
\includegraphics[width=0.49\textwidth]{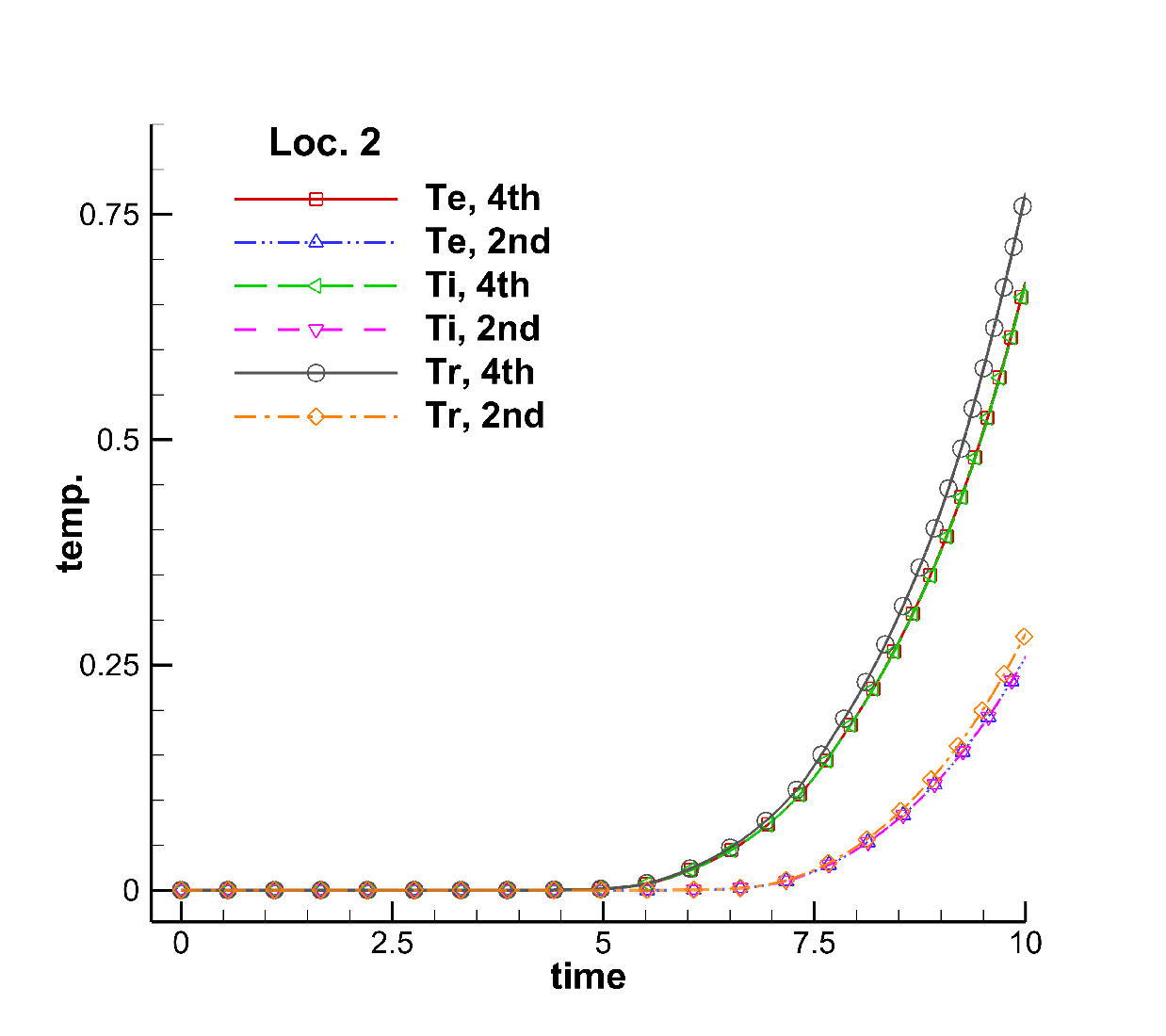}
\includegraphics[width=0.49\textwidth]{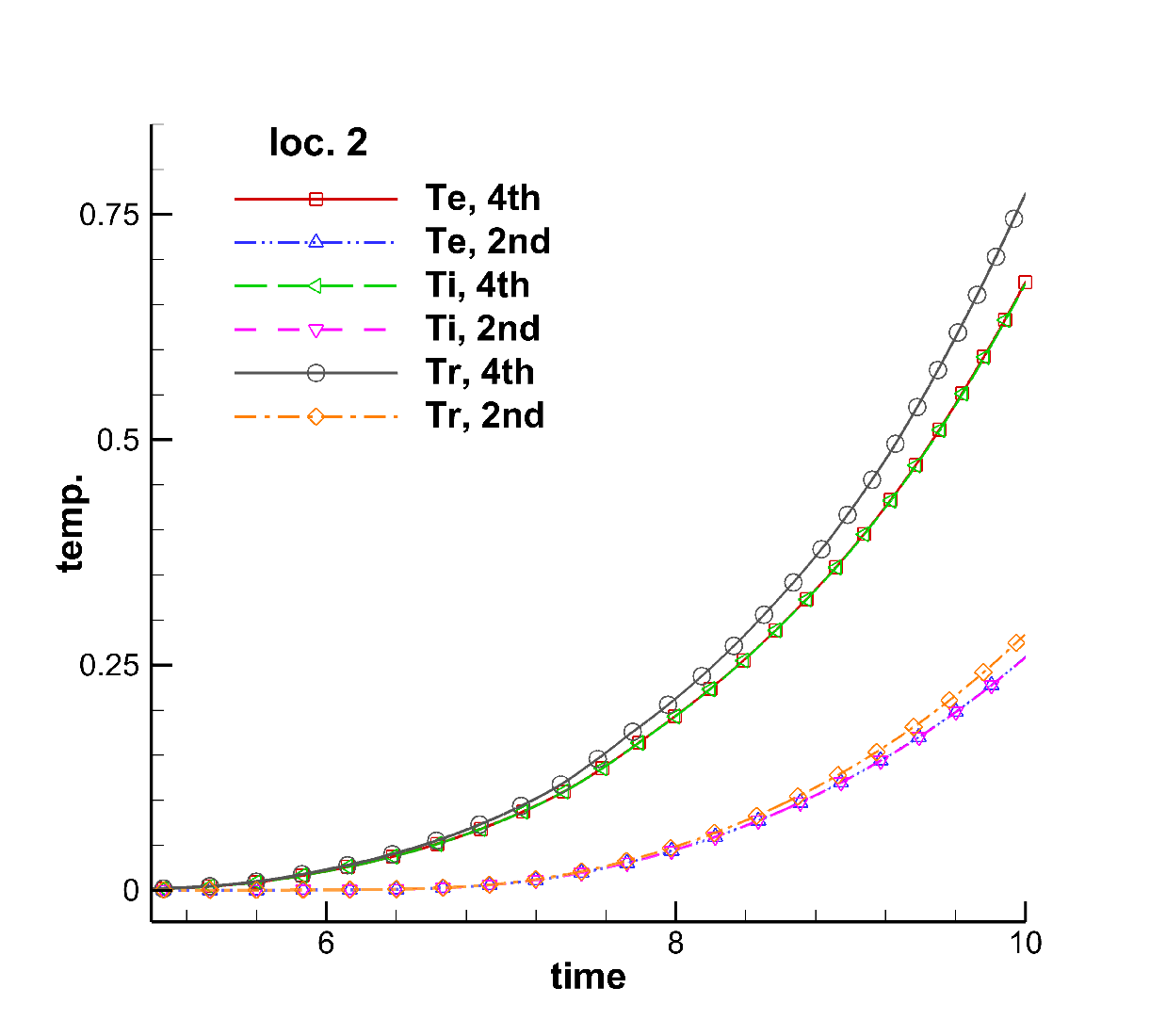}
\caption{\label{3d-ICF-model-2} 3D ICF model problem: Temporal evolution of the three temperatures at Location~1 (top) and Location~2 (bottom), with locally magnified views shown in the right panels. The results compare the performance of the 4th-order central GENO scheme and the linear 2nd-order central scheme. }
\end{figure}

Figure \ref{3d-ICF-model-2} illustrates the temporal evolution of the temperatures at two observation points, computed by the fourth-order central GENO scheme and the linear second-order central scheme using identical time steps. Locally enlarged views are also provided for detailed comparison.
The cell centers corresponding to Locations~1 and~2 are situated at $(0,\,112.5,\,112.5)$ and $(0,\,87.5,\,87.5)$, respectively.
Location~1 is located near the computational domain boundary, which features an initial temperature discontinuity. In contrast, Location~2 is situated in subregion $\Omega_2$ (the middle layer), where material property discontinuities exist.
Due to the extremely large temperature gradient at the isothermal boundary at the initial time, the spatial accuracy of the numerical scheme significantly impacts both the heat flux entering from the boundary and the subsequent temperature rise within the domain.
Compared to the second-order scheme, the fourth-order scheme captures a more accurate, and consequently steeper, temperature gradient.
This leads to a faster temperature increase, as shown in the upper panels of Figure~\ref{3d-ICF-model-2}, and a more rapid propagation of the temperature front, as shown in the lower panels of Figure~\ref{3d-ICF-model-2}.
Furthermore, as the radiation temperature rises, the radiation diffusion coefficient increases dramatically, becoming two orders of magnitude larger than the electron thermal conductivity in subregion $\Omega_3$ adjacent to the isothermal boundary. Consequently, the radiation temperature front propagates significantly faster than the electron and ion temperature fronts. This behavior is consistently captured by both the fourth- and second-order schemes, as illustrated in the bottom panels of Figure \ref{3d-ICF-model-2}.

\begin{table}[htbp]
    \centering
    \begin{tabular}{c|c|c|c|c|c|c}
        \hline
        Reconstruction       & $T_{e,min}$ & $T_{e,max}$ & $T_{i,min}$ & $T_{i,max}$ & $T_{r,min}$ & $T_{r,max}$ \\ \hline
        $t=0.3$             & $3.0\times 10^{-4}$ & 0.522549 & $3.0\times 10^{-4}$ & 0.513476 & $3.0\times 10^{-4}$ & 1.99242 \\ \hline
        $t=5.0$             & $3.0\times 10^{-4}$ & 1.819610 & $3.0\times 10^{-4}$ & 1.819000 & $3.0\times 10^{-4}$ & 1.99940 \\ \hline
    \end{tabular}
    \caption{3D ICF model problem: Temperature bounds at $t=0.3$ and $t=5.0$ computed using the 4th-order GENO scheme with implicit time integration. The results demonstrate that the scheme successfully preserves the physical bounds.}
    \label{3d-ICF-model-table}
\end{table}

\section{Conclusion}

This study presents a high-order finite-volume GENO scheme for the 3TRD equations, specifically designed to handle extreme temperature gradients, solution discontinuities, and stiff diffusion and source terms. The core innovation of this approach lies in the central GENO reconstruction method, which utilizes symmetric sub-stencils.
This method ensures numerical consistency with the isotropic nature of thermal diffusion while achieving the ENO property and facilitating flux evaluation at cell interface, even at material interfaces characterized by discontinuous physical properties.

Furthermore, this study presents the application of a dual time-stepping implicit approach to 3TRD systems. By adopting the implicit temporal method decoupled from complex nonlinear spatial discretization, this strategy overcomes the severe time-step restrictions imposed by stiff operators, thereby enabling large time-step integration with significantly enhanced computational efficiency while preserving spatial high-order accuracy.
Comprehensive numerical validation demonstrates that the proposed scheme achieves the designed high-order accuracy and preserves essential physical bounds, remaining robust even for challenging cases involving temperature discontinuities and time steps orders of magnitude larger than those permitted by explicit schemes. These attributes make the scheme particularly well-suited for demanding applications in high-energy-density physics and inertial confinement fusion simulations.

Future work will extend this methodology to unstructured meshes to accommodate the complex geometries encountered in realistic applications. This effort will require the development of a novel symmetric reconstruction method to ensure the accurate and physics-consistent evaluation of diffusion fluxes, while maintaining algorithmic simplicity and computational efficiency. Additionally, coupling the 3TRD system with hydrodynamics will pave the way for a comprehensive multi-physics simulation framework tailored for advanced scientific computing.

\section*{Acknowledgments}
The authors gratefully acknowledge Dr. Yue Zhang and PhD candidate Hongyu Liu for helpful discussions on the implicit dual time-stepping method.
The current research is supported by National Key R\&D Program of China (Grant Nos. 2022YFA1004500), National Science Foundation of China (92371107), and Hong Kong research grant council (16301222, 16208324).

\section*{References}
\bibliographystyle{ieeetr}
\bibliography{GENO_3TRD}

\section*{Appendix}

This section presents the Jacobian matrix in Eq. (\ref{implicit-3}) for the dual time-stepping implicit time integration of the 3TRD system. The Jacobian matrix is defined as
\begin{align*}
\mathbf{J}=\frac{\mathrm{d} \mathcal{L}}{\mathrm{d} \mathbf{W}},
\end{align*}
where $\mathcal{L}$ is given in Eq. (\ref{3TRD-L}).

The Jacobian matrix of the thermal diffusion term in $\mathcal{L}$ is expressed in component form as
\begin{align*}
\frac{\partial \mathcal{L}^{F}_{\alpha}}{\partial W_{\beta}}=\frac{\mathrm{d} \big(\sum_m \delta F_{\alpha,m}\big)}{\mathrm{d} W_{\beta}},
\end{align*}
and
\begin{align*}
\delta F_{\alpha,m}=\sum_{s=1}^{4}d_{s}\big[\big(k_{\alpha}\frac{\partial T_{\alpha}}{\partial x_m}\big)_{m+1/2,s}- \big(k_{\alpha}\frac{\partial T_{\alpha}}{\partial x_m}\big)_{m-1/2,s} \big]/\Delta x_m,
\end{align*}
where $m\in \{j,k,l\}$ represents the three coordinate directions, $\alpha$ and $\beta$ denote the indices of the three components of vectors $\mathcal{L}^{F}$ and $ \mathbf{W}$, i.e., $\alpha,~\beta \in\{e,i,r\}$, and $d_{s}=1/4$.
Since $\Delta \mathbf{W}$ in Eq. (\ref{implicit-3}) converges to zero in the dual time-stepping method, the discretization of fluxes in $\mathcal{L}^{F}_{\alpha}$ is flexible. Here, we adopt the following second-order approximation
\begin{align*}
\delta F_{\alpha,m}=\big[k_{\alpha,m+1/2}\frac{T_{\alpha,m+1}-T_{\alpha,m}}{\Delta x_m}- k_{\alpha,m-1/2}\frac{T_{\alpha,m}-T_{\alpha,m-1}}{\Delta x_m} \big]/\Delta x_m.
\end{align*}
Then we obtain
\begin{align*}
\begin{split}
&\frac{\partial \mathcal{L}^{F}_{\alpha}}{\partial W_{\beta,m-1}}=\frac{k_{\alpha,m-1/2}}{\Delta x_m^2}\frac{\partial T_{\alpha,m-1}}{\partial W_{\beta,m-1}},\\
&\frac{\partial \mathcal{L}^{F}_{\alpha}}{\partial W_{\beta,m+1}}=\frac{k_{\alpha,m+1/2}}{\Delta x_m^2}\frac{\partial T_{\alpha,m+1}}{\partial W_{\beta,m+1}},\\
&\frac{\partial \mathcal{L}^{F}_{\alpha}}{\partial W_{\beta,m}}  =-\big(\frac{k_{\alpha,m-1/2}+k_{\alpha,m+1/2}}{\Delta x_m^2}\big)\frac{\partial T_{\alpha,m}}{\partial W_{\beta,m}}.
\end{split}
\end{align*}
For the linear model where $W_{\alpha}=c_{\alpha}T_{\alpha}$ and $k_{\alpha}$ is constant, the partial derivative is given by $\partial T_{\alpha}/\partial W_{\beta}=\delta_{\alpha\beta}/c_{\alpha}$.
For the actual physical model, the third component of $\mathbf{W}$ is defined as $W_r=c_r T_r^4$. Thus, the only difference lies in the derivative $\partial T_r/\partial W_r = 1/(4c_r T_r^3)$. Furthermore, although $k_{\alpha}$ is a temperature-dependent function, the terms related to the derivative of $k_{\alpha}$ with respect to $W_{\alpha}$ in the Jacobian matrix are neglected by linearizing the flux.

The Jacobian matrix for the source term in $\mathcal{L}$ is given as
\[
\frac{\mathrm{d} \mathcal{L}^{S}}{\mathrm{d} \mathbf{W}} = \begin{pmatrix}
-(\omega_i+\omega_r)\frac{\mathrm{d} T_e}{\mathrm{d} W_e} & \omega_i\frac{\mathrm{d} T_i}{\mathrm{d} W_i} & \omega_r\frac{\mathrm{d} T_r}{\mathrm{d} W_r} \\[0.5em]
\omega_i\frac{\mathrm{d} T_e}{\mathrm{d} W_e} & -\omega_i \frac{\mathrm{d} T_i}{\mathrm{d} W_i}  & 0 \\[0.5em]
\omega_r\frac{\mathrm{d} T_e}{\mathrm{d} W_e} & 0 & -\omega_r \frac{\mathrm{d} T_r}{\mathrm{d} W_r}
\end{pmatrix}.
\]
In the linear model problem, $\omega_{\alpha}$ is treated as a constant. For the actual physical model, although $\omega_{\alpha}$ is a temperature-dependent coefficient, the terms related to the derivative of $\omega_{\alpha}$ with respect to $W_{\beta}$ in the Jacobian matrix are neglected by linearizing the source term.
Consequently, the Jacobian matrix is determined as $J_{\alpha\beta}=\partial \mathcal{L}^{F}_{\alpha}/\partial W_{\beta} + \mathcal{L}^{S}_{\alpha}/\partial W_{\beta} $.

\end{document}